\documentclass{amsart}
\usepackage{amssymb,latexsym,amsmath}
\usepackage[mathscr]{eucal}
\usepackage{verbatim}
\usepackage{color}
\usepackage{amsthm}
\usepackage{cite}
\usepackage{enumitem}
\theoremstyle{remark}

\def\sideremark#1{\ifvmode\leavevmode\fi\vadjust{\vbox to0pt{\vss
\hbox to 0pt{\hskip\hsize\hskip1em
\vbox{\hsize2cm\tiny\raggedright\pretolerance10000
\noindent#1\hfill}\hss}\vbox to8pt{\vfil}\vss}}}

\newtheorem*{theorem*}{Theorem}
\newtheorem*{corollary*}{Corollary}
\newtheorem{theorem}{Theorem}[section]
\newtheorem{corollary}[theorem]{Corollary}
\newtheorem{lemma}[theorem]{Lemma}
\newtheorem{proposition}[theorem]{Proposition}
\newtheorem{remark}[theorem]{Remark}
\newtheorem{definition}[theorem]{Definition}
\newtheorem{example}[theorem]{Example}
\newtheorem{question}[theorem]{Question}
\newcommand{\SC}{\mathcal{S}}

\newcommand{\be}{\begin{equation}\label} 
\newcommand{\ee}{\end{equation}}
\newcommand{\bq}{\begin{equation*}}
\newcommand{\eq}{\end{equation*}}
\newcommand{\ba}{\begin{align*}}
\newcommand{\ea}{\end{align*}}
\newcommand{\bp}{\begin{proof}}
\newcommand{\ep}{\end{proof}}
\newcommand{\bL}{\begin{lemma}\label}
\newcommand{\eL}{\end{lemma}}
\newcommand{\bP}{\begin{proposition}\label}
\newcommand{\eP}{\end{proposition}}
\newcommand{\bC}{\begin{corollary}\label}
\newcommand{\eC}{\end{corollary}}
\newcommand{\bT}{\begin{theorem}\label}
\newcommand{\eT}{\end{theorem}}
\newcommand{\bTT}{\begin{theorem*}\label}
\newcommand{\eTT}{\end{theorem*}}
\newcommand{\bR}{\begin{remark}\label}
\newcommand{\eR}{\end{remark}}
\newcommand{\bD}{\begin{definition}\label}
\newcommand{\eD}{\end{definition}}

\newcommand{\bE}{\begin{example}\label}
\newcommand{\eE}{\end{example}}

\author{
Sasmita Patnaik\
\and
Gary Weiss\
}

\begin{document}

\title{Interplay of Simple and Selfadjoint-ideal semigroups in $B(H)$}
\maketitle

\begin{abstract}
This paper investigates a question of Radjavi: Which multiplicative
semigroups in $B(H)$ have all their ideals selfadjoint (called herein
selfadjoint-ideal (SI) semigroups)? We proved this property is a unitary invariant for $B(H)$-semigroups, which invariant we believe is new. 

We characterize those SI semigroups $\SC$ singly generated by $T$,  for $T$ a normal operator and for $T$ a rank one operator. When $T$ is nonselfadjoint and normal or
rank one: $\SC$ is an SI semigroup if and only if it is simple, except in one
special rank one partial isometry case when our characterization yields $\SC$
that are SI but not simple. So SI and simplicity are not equivalent notions. When $T$ is selfadjoint, it is straightforward to see that $\SC$ is always an SI semigroup, but we prove by examples they may or may not be simple, but for this case we do not have a characterization. 

The study of SI semigroups involves solving certain operator equations in the
semigroups. A central theme of this paper is to study when and when not SI is
equivalent to simple.

\end{abstract}

\section{Introduction}
We investigate the recent question of Radjavi: Which multiplicative semigroups  in $B(H)$ have all their multiplicative ideals (that is, semigroup-ideals) automatically selfadjoint (Definitions \ref{D1}-\ref{D4}).  Like simplicity, this turns out to be a unitary invariant property for $B(H)$-semigroups (Theorem \ref{T1}), which invariant we believe is new. We will call these semigroups selfadjoint-ideal semigroups (SI semigroups for short).
Herein, all semigroup-ideals are two-sided, and $B(H)$ denotes the set of bounded linear operators on a finite or infinite dimensional complex Hilbert space $H$.

In the study of semigroups, it appears to us unexplored the possibility of general characterizations (or any general structure) of selfadjoint semigroups of $B(H)$ (i.e., semigroups closed under adjoints, see Definition \ref{D1}).  This complicates our study of those selfadjoint semigroups whose ideals are all selfadjoint. The latter question on selfadjointness of all ideals in a semigroup $\SC$ is equivalent to solving for each $T\in \SC$,  $T^* = XTY$ for $X, Y \in \SC \cup \{I\}$ for non-unital $\SC$ (or for $X, Y \in \SC$ when $\SC$ is unital) (Lemma \ref{L2}). Radjavi first noticed that this operator equation is solvable for $B(H)$ and $F(H)$ regarded as semigroups. $F(H)$ denotes the finite-rank operators in $B(H)$.

There is continuing interest in the study of  multiplicative semigroups in $B(H)$ from different perspectives. A common theme has been to characterize special classes of $B(H)$-semigroups (that is, semigroups of $B(H)$), e.g., \cite{JLPR, C,O,OR,OMR2,OMR,P,R,RR,RRS,S}. For instance, Radjavi and Popov linked selfadjoint semigroups of partial isometries to a representation of inverse semigroups \cite{RP}. (An inverse semigroup is a semigroup in which every element $x$ has a unique inverse $y$ in the semigroup in the sense that $x = xyx$ and $y = yxy$, \cite{DP}, \cite{P}.)  In a sequel to \cite{RP}, Radjavi, Popov, et al., more recently broadened their investigation by determining sufficient conditions under which selfadjoint extensions of semigroups of partial isometries (i.e., semigroups generated by a semigroup of partial isometries and its adjoint) are again semigroups of partial isometries\cite[Theorem 2.6]{JLPR}.  
 In this paper we study semigroups in $B(H)$ that possess Radjavi's recently introduced  automatic selfadjoint-ideals property mentioned above, and along the way we expand our study of $B(H)$-semigroups (i.e., semigroups of $B(H)$) beyond partial isometry semigroups. 

The first easy examples of selfadjoint $B(H)$-semigroups possessing this selfadjoint-ideal property that we encountered were $F(H)$ (Corollary \ref{C1}), $B(H)$ (by the same argument used for $F(H)$), and selfadjoint semigroups of unitaries (Remark \ref{R3}(iv)). A less trivial example is a singly generated selfadjoint semigroup $\SC(V, V^*)$ generated by an isometry $V$ (Remark \ref{R3}(v)). Other historical examples of characterizations of semigroups: one-parameter semigroups of unitary operators were characterized by Stone \cite{S} as they are related to certain selfadjoint operators, and the one-parameter semigroups of isometries were characterized by Cooper \cite{C} as they are related to certain symmetric operators in Hilbert space. 

%One-parameter semigroups of unitary operators were characterized by Stone \cite{S} as they are related to certain selfadjoint operators, and the one-parameter semigroups of isometries were characterized by Cooper \cite{C} as they are related to certain symmetric operators in Hilbert space. (An unbounded operator $T$ is a symmetric operator if $\left<Tx, y\right> = \left<x, Ty\right>$ for all $x,y$ in the domain of $T$.)

The next natural class of operators to consider  for our study of selfadjoint-ideal semigroups are partial isometries.
These led us to some illuminating examples for the development of a theory. 

Notably, the full class of partial isometries is a selfadjoint class but does not form a semigroup as it is not closed under products. The product $VW$ of two partial isometries $V$ and $W$ is a partial isometry if and only if the projections $V^*V$ and $WW^*$ commute \cite[Lemma 2]{HW}. Two partial isometries in $M_2(\mathbb C)$ whose product is not a partial isometry are 
\begin{equation*}
V = \begin{pmatrix}
0&0\\
1&0
\end{pmatrix}
\quad \text{and}\quad W = \begin{pmatrix}
\frac{1}{\sqrt{2}} & 0\\
\frac{1}{\sqrt{2}} &0
\end{pmatrix}.
\end{equation*}
The matrices $V$ and $W$ are partial isometries, but their product 
\begin{equation*}
VW = \begin{pmatrix}
0&0\\
\frac{1}{\sqrt{2}}&0
\end{pmatrix}
\end{equation*}
 is not a partial isometry because $VW \neq VW(VW)^{*}VW$, recalling $U$ is a partial isometry if and only if $U = UU^{*}U$ \cite[Corollary 3 of Problem 127]{H}. More easily but less significant, $VW$ is not a nonzero partial isometry because $\frac{1}{\sqrt{2}} = ||VW|| \neq 1$. 
 
Historically,  it appears to us that Halmos and Wallen initiated the study of partial isometries in their seminal paper \textit{Powers of Partial Isometries} \cite{HW}  where they completely described the structure of power partial isometries (those partial isometries $T \in B(H)$ for which all $T^{n}$ are also partial isometries): Every power partial isometry is a direct sum whose summands are unitary operators, pure isometries, pure co-isometries, and truncated shifts \cite[Theorem]{HW}. (They understood that not all the four summands must be present in each case allowing the possibility of $0$ as one of the summands [ibid, p. $657$ last paragraph].) Halmos and Wallen state that it is not immediately obvious that for $T$ a partial isometry, $T^{2}$ may not be a partial isometry \cite[Example]{HW}. 
%As motivation for [ibid, Theorem] characterizing the structure of power partial isometries  in their introduction, they commented ``a unitary operator $U$ has a skew-Hermitian logarithm ($U = \exp(iH)$, $H$ selfadjoint), and an isometry is a direct sum of a unitary operator and a unilateral shift with a suitable multiplicity (this is commonly called the Wold decomposition).''

Radjavi and Popov observed the following connection between a power partial isometry $T$ and the singly generated selfadjoint semigroup $\SC(T, T^*)$ generated by $T$: 
%(for this notation see last line of the paragraph following Definition \ref{D6}, also \cite[second last paragraph of the Introduction]{RP}): 
An operator $T$ is a power partial isometry if and only if the semigroup $\SC(T, T^*)$ consists of only partial isometries \cite[Proposition 2.2]{RP}. We observed that this result also connects power partial isometries $T$ with the selfadjoint property of the ideals of $\SC(T, T^*)$, because whenever $\SC(T, T^*)$ consists of only partial isometries, $\SC(T, T^*)$ is a selfadjoint-ideal semigroup (Corollary \ref{CPP}). So we have, if $T$ is a power partial isometry, then $\SC(T, T^*)$ is a selfadjoint-ideal semigroup. In this case, we shall see later that $\SC(T, T^*)$ may not be simple (Example \ref{ENS}). However, under special conditions (e.g. rank one power partial isometries), $\SC(T, T^*)$ is simple (Corollary \ref{CPI}). But for a particular subclass of power partial isometries, the selfadjoint-ideal semigroup $\SC(T, T^*)$ is not simple (Example \ref{ENS}). A connection we found between the notion of SI semigroup and power partial isometry is: For $T$ a normal nonselfadjoint operator, $\SC(T, T^*)$ is a selfadjoint-ideal semigroup if and only if $T$ is a power partial isometry of the form $U \oplus 0$ ($U$ unitary) if and only if $\SC(T, T^*)$ is simple (Theorem \ref{N2}). 

In this investigation it is natural to first study general singly generated selfadjoint semigroups $\SC(T, T^*)$.
The focus of this paper is to characterize when they are selfadjoint-ideal semigroups for $T$ a normal operator and for $T$ a rank one operator.  Necessary and sufficient conditions for these semigroups to be selfadjoint-ideal semigroups are given in Theorem \ref{N2}, Theorems \ref{TR1}-\ref{TN1}, and Theorem \ref{TR2}. Unexpected results on simplicity from this investigation are that the class of selfadjoint-ideal semigroups $\SC(T, T^*)$ generated by rank one partial isometries $T$ with $0 < tr T <1$ are not simple; but the class of selfadjoint-ideal semigroups $\SC(T, T^*)$ generated by a normal operator, or by a rank one partial isometry with $|tr T|=1$, or by a rank one non-partial isometry are simple.
Upshot: singly generated selfadjoint-ideal semigroups $\SC(T, T^*)$ generated by normal and rank one operators $T$ are usually simple and rarely nonsimple. The selfadjoint-ideal property versus simplicity for higher rank operators remains to be studied.

A summary at the end of Section 3 is given of the main characterizations in Sections 2-3 leading to a  classification of $\SC(T, T^*)$ for simple and nonsimple, for $T$ rank one and $T$ normal  with one open question exception (selfadjoint with rank greater than one).
\newpage

\textbf{Terminology (Definitions \ref{D1}-\ref{D6})}
\vspace{.2cm}

The terminology given in Definitions \ref{D1}-\ref{D3} is standard (as in \cite{SC},  \cite{RP}, and \cite{RR}). The terminology in Definition \ref{D4} on the notion of selfadjoint-ideal semigroups and our focus in this paper, we believe is new and due to Radjavi.

\bD{D1}
A semigroup $\SC$ in $M_n(\mathbb C)$ or $B(H)$ is a subset closed under multiplication. \linebreak  A selfadjoint semigroup $\SC$ is a semigroup also closed under adjoints,  i.e., $\SC^* := \{T^* \mid T \in \SC\} \subset \SC$. 
\eD
 
\bD{D3}
 An ideal $J$ of a semigroup $\mathcal S$ in $M_{n}(\mathbb C)$ or $B(H)$ is a subset of $\SC$ closed under products of operators in $\SC$ and $J$. That is, $XT, TY \in J$ for $T \in J$ and $X, Y \in \mathcal S$. And so also $XTY \in J$. 
 \eD

   \bD{D4} A selfadjoint-ideal (SI) semigroup $\SC$ in $M_{n}(\mathbb C)$ or $B(H)$ is a semigroup for which every ideal $J$ of $\SC$ is closed under adjoints, i.e., $J^{*} := \{T^{*} \mid T \in J\} \subset J.$  
 \eD
Because this selfadjoint ideal property in Definition \ref{D4} concerns selfadjointness of all ideals in a semigroup, we call these semigroups  selfadjoint-ideal semigroups (SI semigroups for short). By Definition \ref{D4} every SI semigroup is itself a selfadjoint semigroup, so SI semigroups is a subclass of selfadjoint semigroups. And it is a proper subclass because Example \ref{EE} below provides easy examples of selfadjoint semigroups that are not SI semigroups. The most obvious category of SI semigroups are semigroups consisting of only selfadjoint operators, e.g., the  semigroup $\SC(T)$ generated by a selfadjoint operator $T$ (that is, $\SC(T) = \{T, T^2, \cdots\}$). They are SI semigroups because they consist of only selfadjoint operators so all their nonzero ideals are selfadjoint. (Note $\SC(T)$ is formally introduced in paragraph following Definition \ref{D5} below.) 

We defined SI semigroup as selfadjointness of all its ideals. So a vacuous case of this are simple selfadjoint semigroups, i.e., selfadjoint semigroups without proper nonzero ideals. However, our ``simple semigroups" may or may not contain zero, unlike standard usage which differentiates between simple and  $0$-simple semigroups. We do this because our focus is the study of semigroups with the SI property which is not affected by the  presence or absence of zero since $\{0\}$ is always a selfadjoint ideal whenever $0$ belongs to the semigroup.
\textit{So in particular a simple selfadjoint semigroup $\SC$ is always an SI semigroup because the semigroup itself is its only nonzero ideal and when $0 \in \SC$, $\{0\}$ is always a selfadjoint ideal.}

Remark \ref{R3}(i) provides examples of SI semigroups which are not simple. In Sections 2-3, we investigate, in several instances, conditions under which singly generated selfadjoint semigroups are simple (\color{blue}\color{black} Theorem \ref{T2}, Theorem \ref{TR1}, and Theorem \ref{TR2}). Simplicity (i.e., no proper ideals) is clearly a stronger condition than SI semigroupness. In summary,
$$\text{simple selfadjoint semigroups $\subsetneq$ SI semigroups $\subsetneq$ selfadjoint semigroups}$$

\vspace{.3cm}

\textit{Semigroups generated by $\mathcal A \subset B(H)$  \cite[Introduction]{RP}}

\bD{D5}
The semigroup generated by a set $\mathcal A \subset  M_n(\mathbb C)$ or $B(H)$, denoted by $\mathcal S(\mathcal A)$, is the intersection of all semigroups containing $\mathcal A.$ Also define $\mathcal A^*:= \{A^* | A \in \mathcal A\}$.
\eD
For short we denote by $\SC(T)$ the semigroup generated by $\{T\}$ (called generated by $T$ for short).
It should be clear for the semigroup $\mathcal S(\mathcal A)$ that Definition \ref{D5} is equivalent to the semigroup consisting of all possible words of the form $A_1A_2\cdots A_k$ where $k \in \mathbb N$ and $A_i \in \mathcal A$ for each $1\leq i \leq k$. 
 
 \bD{D6}
The selfadjoint semigroup generated by a set $\mathcal A \subset M_n(\mathbb C)$ or $B(H)$ denoted by $\SC(\mathcal A \cup \mathcal A^{*})$ or $\SC(\mathcal A,\mathcal A^{*})$, is the intersection of all selfadjoint semigroups containing  $\mathcal A \cup \mathcal A^{*}$. 

Let $\SC(T, T^*)$ denote for short $\SC(\{T\}, \{T^*\})$ and call it the singly generated selfadjoint semigroup generated by $T$.
 \eD
 
It is clear that $\SC(\mathcal A, \mathcal A^{*})$ is a selfadjoint semigroup. Moreover, it is clear that Definition \ref{D6} conforms to the meaning of $\SC(\mathcal A \cup \mathcal A^*)$ in terms of words discussed above. That is, it consists of all words of the form $A_1A_2\cdots A_k$ where $k \in \mathbb N$ and $A_i \in \mathcal A \cup \mathcal A^{*}$ for each $1\leq i \leq k$. \\

The connection between $\mathcal S(T)$ and $\SC(T, T^*)$ in relation to SI is given later in Proposition \ref{PS}.\\

Next in Proposition \ref{P1} below we provide an explicit list for the operators in the semigroup $S(T, T^*)$ on which we  depend in two key places of our work (proof of Theorem \ref{T2} and Remark \ref{RS}).\\

As stated in the Introduction (first paragraph), throughout $H$ can be finite or infinite dimensional unless otherwise specified. Unless relevant, we will not mention separately cases $M_n(\mathbb C)$ and $B(H)$.

\bP{P1}For $T \in B(H)$, the semigroup $S(T, T^*)$ generated by the set $\{T, T^*\}$ is given by 

\noindent $ S(T, T^*) = \{T^n, {T^*}^n, \Pi_{j=1}^{k}T^{n_j}{T^*}^{m_j},  (\Pi_{j=1}^{k}T^{n_j}{T^*}^{m_j})T^{n_{k+1}}, \Pi_{j=1}^{k}{T^*}^{m_j}T^{n_j}, (\Pi_{j=1}^{k}{T^*}^{m_j}T^{n_j}){T^*}^{m_{k+1}},$ where $n \ge 1,\,  k\ge1,\, n_j, m_j \ge 1\, \text{for}\, 1 \le j \le k, ~\text{and}~n_{k+1}, m_{k+1} \geq 1\}$.
The product $\Pi_{j= 1}^{k}$ in the semigroup list is meant to denote an ordered product. 
\eP
\bp
This follows directly from the paragraph on words mentioned after Definition \ref{D6} by taking $\mathcal A = \{T\}$.
\ep

Alternatively $\SC(T, T^*)$ consists of: words only in $T$, words only in $T^*$, words that begin and end in $T$, words that begin with $T$ and end with $T^*$, and words that begin with $T^*$ and end with $T$ and words that begin and end with $T^*$.\\ 
 
 \textit{Ideals generated by $\mathcal A \subset \SC$ for $\SC$ a semigroup in $B(H)$} 
 \vspace{.3cm}
 
Recall that we consider only two-sided semigroup-ideals (first paragraph of the Introduction), that is, a subset of a semigroup closed under multiplication between elements of the subset and the elements of the semigroup. As defined in ring theory and semigroup theory \cite{G}, we define an ideal of a semigroup $\SC$ generated by a subset $\mathcal A \subset \SC$, denoted as $(\mathcal A)_{\SC}$, to be the intersection of all the ideals of $\SC$ containing $\mathcal A$. 
It is easy to verify that arbitrary intersections of ideals of a semigroup is again an ideal of that semigroup. Therefore $(\mathcal A)_{\SC}$ is the smallest ideal of $\SC$ containing $\mathcal A$. 
We call $(\mathcal A)_{\SC}$ the ideal of $\SC$ generated by a subset $\mathcal A$.
One has a simple characterization of $(\mathcal A)_{\SC}$  as follows.

\bL{L1}
For $\SC$ a semigroup and $\mathcal A \subset \mathcal S$, $$(\mathcal A)_{\mathcal S} = \mathcal A \cup \mathcal S\mathcal A \cup \mathcal A\mathcal S \cup \SC \mathcal A\SC$$ where $\mathcal S\mathcal A := \{XA | X \in \SC, A \in \mathcal A\}$, $\mathcal A\mathcal S := \{AX | X \in \SC, A \in \mathcal A\}$, $\SC \mathcal A \SC := \{XAY | X, Y \in \SC, A \in \mathcal A\}$. 
 \eL
\bp{}
Since $(\mathcal A)_{\mathcal S}$ is an ideal of $\mathcal S$ containing $\mathcal A$, by Definition \ref{D4}, $\mathcal A \cup \mathcal S\mathcal A \cup \mathcal A\mathcal S \cup \SC \mathcal A \SC \subset (\mathcal A)_{\mathcal S}$. Also $\mathcal A \cup \mathcal S\mathcal A \cup \mathcal A\mathcal S \cup \SC \mathcal A \SC$ is clearly an ideal of $\mathcal S$ containing $\mathcal A$. So $(\mathcal A)_{\mathcal S}$ being the smallest ideal containing $\mathcal A$, one has $(\mathcal A)_{\mathcal S} \subset \mathcal A \cup \mathcal S\mathcal A \cup \mathcal A\mathcal S \cup \SC \mathcal A \SC$, hence one obtains equality.
\ep

\bL{LP} {(An equivalent condition for the selfadjointness of an ideal generated by a set)}\\
For $\SC$ a selfadjoint semigroup and $\mathcal A \subset \mathcal S$, the ideal $(\mathcal A)_{\SC}$ is selfadjoint if and only if $A^* \in (\mathcal A)_{\SC}$ for every $A \in \mathcal A$.
 \eL
 \bp
 Clearly, if the ideal $(\mathcal A)_{\SC}$ is selfadjoint, then $A^* \in (\mathcal A)_{\SC}$ for every $A \in \mathcal A$. Conversely, by Lemma \ref{L1}, $$(\mathcal A)_{\mathcal S} = \mathcal A \cup \mathcal S\mathcal A \cup \mathcal A\mathcal S \cup \SC \mathcal A\SC.$$
 So an operator in the ideal  $(\mathcal A)_{\SC}$ is of the form $XAY$ where $X, Y \in \SC \cup \{I\}$ and $A \in \mathcal A$. Therefore $(XAY)^* = Y^*A^*X^* \in (\mathcal A)_{\SC}$ as $A^* \in (\mathcal A)_{\SC}$ and $X^*, Y^* \in \SC$.    \ep

\bL{L2}
For semigroup $\SC \subset B(H)$, the following are equivalent:

\begin{enumerate}[label=(\roman*)]
\item $\SC$ is an SI semigroup. 
\item Every principal ideal of $\SC$ is selfadjoint. In symbols, for all $T \in \SC$, $T^* \in (T)_{\SC}$.
\item For all $T \in \SC$, $T^* = S_1TS_2$ for some $S_1, S_2 \in \SC \cup \{I\}$. 
\end{enumerate}
\eL
\bp{}
$\text{(i)} \Rightarrow \text{(ii)}$
Follows directly from the definition of an SI semigroup (Definition \ref{D4}). \\

(ii) $\Rightarrow$ (iii) Applying Lemma \ref{L1} to $\mathcal A = \{T\}$, $$(T)_{\SC} = \{T\} \cup  \mathcal S\{T\} \cup \{T\}\mathcal S \cup \SC \{T\} \SC.$$ Therefore $T^* \in (T)_{\SC}$ or equivalently $T^* = S_1TS_2$ for some $S_1, S_2 \in \SC \cup \{I\}$.\\

$\text{(iii)} \Rightarrow \text{(i)}$ It suffices to show that every nonzero ideal $J$ is selfadjoint (as $\{0\}$ is always selfadjoint, if it happens that $0 \in \SC$). For $0 \neq T \in J$, $(T)_{\SC} \subset J$ as $J$ is an ideal. If $T$ is selfadjoint, then $T^* = T \in J$. Assume $T$ is nonselfadjoint. For each $T \in \SC$, by (iii) one has $T^* = S_1TS_2$ where at least  $S_1$ or $S_2 \in \SC \setminus \{I\}$ (since $T \neq T^*$). Since $J$ is an ideal containing $T$, $S_1TS_2 \in J$ so $T^* \in J$. Hence $J$ is selfadjoint. Therefore $\SC$ is an SI semigroup.
\ep

In consideration of Lemma \ref{L2}, it is natural to ask if there is a connection between a semigroup $\SC$ possessing the SI property and all the singly generated semigroups $\SC(T, T^*)$ of its members $T$ possessing SI. The next Proposition \ref{PS} shows the connection.  However, admittedly we yet have no example of a semigroup with this property. 
\bP{PS} (A sufficient but not a necessary condition for $\SC$ to be an SI semigroup.)\\
$\SC$ is an SI semigroup if all its singly generated selfadjoint semigroups $\SC(T, T^*)$ are SI semigroups. 
\eP
\bp
For $T \in \SC$ and principal ideal $(T)_{\SC}$, one has $\SC(T, T^*) \subset \SC$ and $(T)_{\SC(T, T^*)} \subset (T)_{\SC} \subset \SC$. Since $\SC(T, T^*)$ is an SI semigroup, $(T)_{\SC(T, T^*)}$ is a selfadjoint ideal. And so by Lemma \ref{L2} (ii) $\Rightarrow$ (iii) applied to $\SC(T, T^*)$, one has $T^* = XTY$ for some $X, Y \in \SC(T, T^*) \cup \{I\}$, so $T^* \in (T)_{\SC}$.  
By Lemma \ref{L2} (ii)$\Rightarrow$(i), $\SC$ is an SI semigroup.

To see the sufficient condition is not necessary, the semigroup $\SC = M_{2}(\mathbb C)$ being selfadjoint and $\SC(T, T^*)$ of Example \ref{E} provides a counterexample to the converse.
\ep

It is also natural to consider the relationship between semigroups and their principal ideals in determining SI instead of having to consider all its ideals. This is given in Remark \ref{sa}(ii) below.
\bR{sa} 
(i) In a selfadjoint semigroup, ideals generated by selfadjoint operators are selfadjoint. Indeed, for $\SC$ a selfadjoint semigroup, let $\mathcal A \subset \SC$ be a subset of selfadjoint operators, hence $\mathcal A = \mathcal A^*$. Then the selfadjointness of the ideal $(\mathcal A)_{\SC}$ follows directly from Lemma \ref{LP}. 

(ii) For a selfadjoint semigroup $\SC$ to be an SI semigroup, it suffices to show that every principal ideal generated by a nonselfadjoint operator is selfadjoint. (Follows immediately from Lemma \ref{L2} (i) $\Leftrightarrow$ (ii) and Remark \ref{sa}(i)).

Moreover, for $ S \in \SC(T, T^*)$, the principal ideal $(S)_{\SC(T, T^*)}$ is selfadjoint if and only if $S^* \in (S)_{\SC(T, T^*)}$ (equivalently, $S^* = XSY$ for some $X, Y \in \SC(T, T^*) \cup \{I\})$. Indeed, the ``if and only if" part follows from Lemma \ref{LP} for $\mathcal A = \{S\}$, and the parenthetical equivalently part follows from Lemma \ref{L1} for $\mathcal A = \{S\}$. But $S^* = XSY \Leftrightarrow S = Y^*S^*X^*$. From this and Lemma \ref{LP} it follows that $(S)_{\SC(T, T^*)}$ is selfadjoint  if and only if $(S^*)_{\SC(T,T^*)}$ is selfadjoint if and only if $(S)_{\SC(T, T^*)} = (S^*)_{\SC(T,T^*)}$.

\eR

Another natural consideration is the relationship between single generated semigroups $\SC(T)$ and  the singly generated selfadjoint semigroup $\SC(T, T^*)$ generated by $T$. For this we have:
\bP{Psa}
For the singly generated semigroup $\SC(T)$ and the singly generated selfadjoint semigroup $\SC(T, T^*)$ generated by $T$, $$\SC(T) ~ \text{is an SI semigroup if and only if} ~\SC(T) = \SC(T, T^*) ~ \text{and}~ \SC(T, T^*)~ \text{is an SI semigroup}.$$
\eP
\bp
Sufficiency is clear. For the necessity, $\SC(T)$ is an SI semigroup then observe that $T^* \in (T)_{\SC(T)} \subset \SC(T) \subset \SC(T, T^*)$. Therefore also $\SC(T, T^*) \subset \SC(T)$ hence equality and which is hypothesized to be an SI semigroup.
\ep
\bR{Rsa} (i) For $T$ a selfadjoint operator, $\SC(T, T^*)$ is trivially an SI semigroup as $\SC(T, T^*) = \SC(T) = \{T^n \mid n\geq 1\}$ where every element is selfadjoint so all its principal ideals are automatically selfadjoint (Remark \ref{sa}(i)). \\
(ii) The SI semigroup $\SC(T, T^*) = \{T^n \mid n\geq 1\}$ is isomorphic to some semigroup of nonnegative integers. For instance, if $T$ is a projection, $T = T^* = T^2$, then  $S(T, T^*) =  \{T\}$ is isomorphic to the semigroup $\{1\}$ under multiplication, and as is well known, in general, infinite singly generated abstract semigroups are isomorphic to the set of all positive integers  under addition.
\eR

In deference to the role partial and power partial isometries played in the work of Radjavi and Popov we have:
\newpage
\bP{PPSI} A semigroup $\SC$ consisting of only partial isometries is an SI semigroup.
\eP
\bp By Lemma \ref{L2} (i) $\Rightarrow$ (ii), if $T$ is in the selfadjoint semigroup $\SC$, it suffices to show the principal ideal it generates, $(T)_{\SC}$, is selfadjoint. Use again the partial isometry characterization, $T = TT^*T$, hence $T^* = T^* TT^* \in (T)_{\SC}$.
\ep
\bC{CPP}
If $T$ is a power partial isometry, then $\SC(T, T^*)$ is SI.
\eC
\bp
Follows by combining Proposition \ref{PPSI} with the Radjavi-Popov result that an operator $T$ is a power partial isometry if and only if the semigroup $\SC(T, T^*)$ consists of only partial isometries \cite[Proposition 2.2]{RP}.
\ep
The converses of Propositions \ref{PPSI} and \ref{PPP} below sometimes fails, that is, there are SI semigroups that do not contain any partial isometry. For instance, consider a selfadjoint operator $T$ with $||T|| <1$, then $\SC(T, T^*) = \{T^n \mid n\geq 1\}$ is clearly SI (Remark \ref{sa}(i)), but none of the elements are partial isometries as the norm of each element is strictly less than one.\\

\textbf{Algebraic $B(H)$-ideals viewed as SI semigroups} \\

An algebraic $B(H)$-ideal, i.e., a two-sided ideal of the ring $B(H)$, can be viewed as a semigroup naturally inheriting the multiplicative structure of the ideal.
We began our investigation of the SI property by considering Radjavi's question on $B(H)$-ideals: Which nonzero proper algebraic $B(H)$-ideals are SI semigroups? A complete answer to this is provided in Corollary \ref{C1}: The finite rank ideal $F(H)$ is the only one.\\

\noindent We need the following lemma in the proof of Theorem \ref{L4}.
\bL{L}\cite[Lemma 1]{FR}
If $K$ and $L$ are positive compact operators, then TFAE:\\
\noindent (i) There exist operators $A$ and $B$, at least  one of them compact, for which $L = AKB.$\\
\noindent (ii)  $s_n(L) = o(s_n(K))$ as $n \rightarrow \infty.$
\eL

\bT{L4}
All semigroups $\SC \subset B(H)$ of compact operators with at least one nonselfadjoint infinite rank compact operator are not SI.
\eT
\bp{}
Suppose otherwise that $\SC$ is an SI semigroup with $T \in \SC$,  a nonselfadjoint compact operator of infinite rank. Then the principal ideal $(T)_{\SC}$ is selfadjoint so $T^{*} \in (T)_{\SC}$. Since $T^{*} \neq T$, by Lemma \ref{L2} (i) $\Rightarrow $(iii) we have $T^{*} = XTY$ for some $X, Y \in \SC \cup {\{I\}}$ and $X$ or $Y \in \SC \setminus \{I\}$. Using the polar decomposition $T = V|T|$, one obtains $|T| = V^{*}T = T^{*}V.$  Replacing $T^{*} = XTY$ in the last product yields $|T| = XTYV = XV|T|YV$. Then $|T| = A|T|B$ for $A = XV, B = YV \in K(H)$ (the ideal of compact operators). Applying Lemma \ref{L} for $L =K = |T|$, one obtains a contradiction to Lemma \ref{L} (ii) because $\forall n, 0< s_n(L) = s_n(|T|) = s_n(K)$ showing $s_n(L) \neq o(s_n(K))$ as $n \rightarrow \infty.$
\ep

\bC{C1} $F(H)$ is an SI semigroup and all other nonzero proper $B(H)$-ideals are not SI semigroups.
\eC
\bp{}
Suppose $J$ is an ideal of the semigroup $F(H)$ and $ T \in J$. As $T = V|T|$ (with ran $V$ $=$ ran $T$ and $V^*T = |T|$) where $V$ is a finite rank partial isometry in the polar decomposition of $T$, one has $T^{*}= V^{*}TV^{*}$ and $V^{*} \in F(H)$. Therefore $T^{*} \in J$, showing $J$ is selfadjoint. Therefore $F(H)$ is an SI semigroup.

Suppose $\SC  \neq F(H)$ is a nonzero proper algebraic $B(H)$-ideal, then it properly contains all finite rank operators \cite[Theorem 1.1, Chapter III]{GK}. Then there is an infinite rank nonselfadjoint compact operator $T \in \SC$. Indeed, for $T \in \SC\setminus F(H)$,  either $T$ or $iT$ is nonselfadjoint and $\SC$ being an algebraic $B(H)$-ideal, also $iT \in \SC$. By Theorem \ref{L4}, $\SC$  is not an SI semigroup.
\ep

In Proposition \ref{PPP}, we next prove that any general $B(H)$-semigroup containing all partial isometries is SI. An easy consequence is that the multiplicative semigroup $B(H)$ is SI (first observed by Radjavi).

\newpage
\bP{PPP} A semigroup $\SC$ containing every partial isometry is an SI semigroup.
\eP
\bp{}
For $T \in \mathcal S$, from the polar decomposition $T = V|T|$ where $V$ is a partial isometry and hence $V \in \SC$, and the fact that also $V^{*}T = |T|$, one has $T^{*} = |T|V^{*} = V^{*}TV^{*}$ and $V^*$ a partial isometry so $V^* \in \SC$. Therefore $T^{*} = V^{*}TV^{*}$. Then for every ideal $J$ in $\SC$, $T^{*} = V^{*}TV^{*} \in J$  for all $T \in J$ since $V^{*} \in \mathcal S$. The fact that when $V$ is a partial isometry so also $V^*$ is a partial isometry comes from the characterization \cite[Corollary 3 of Problem 127]{H}: $T$ is a partial isometry iff $T = TT^*T$.
\ep

For finite dimensional Hilbert space, the ring $M_{n}(\mathbb C)$ is well-known to contain no nonzero algebraic ideal (i.e., ideal in the ring). Contrasting ring ideal structure with semigroup-ideal structure, $M_{n}(\mathbb C)$, $n \ge 2$, regarded as a multiplicative semigroup, contains non-trivial semigroup-ideals. For instance, the principal semigroup-ideal $(T)_{\SC(T, T^*)}$ generated by a nilpotent operator $T \in M_{n}(\mathbb C)$ (so $det T = 0$), is proper because $I \not \in (T)_{\SC(T, T^*)}$ since otherwise $XTY = I$ and hence the contradiction $0 = det X det T det Y = det(XTY) = det I = 1$. Hence $\{0\} \neq (T)_{M_n(\mathbb C)} \neq M_n(\mathbb C)$, i.e.,  $M_n(\mathbb C)$ has non-trivial ideals. It then becomes natural to study multiplicative semigroups of $M_{n}(\mathbb C)$ and their semigroup-ideals in view of the role the SI property plays in their structure. 

For infinite dimensional Hilbert space, the relationship between algebraic $B(H)$-ideals and multiplicative $B(H)$-semigroups is as follows. A  proper two-sided algebraic $B(H)$-ideal is also a selfadjoint semigroup of $B(H)$ with respect to multiplication and each one consists of compact operators only. 
And among these algebraic $B(H)$-ideals, only $F(H)$ is an SI semigroup as shown in Corollary \ref{C1} above. This also answers Radjavi's natural question to determine which proper algebraic $B(H)$-ideals are SI semigroups? 

\bR{R-SI}
In summary, among all algebraic $B(H)$-ideals, $F(H)$ and $B(H)$ are the only semigroups that are SI. Moreover, viewed as multiplicative semigroups, both are not simple. Indeed, the principal ideal generated by a rank one orthogonal projection is a proper semigroup-ideal in the semigroup $F(H)$ because every operator in that ideal is of rank at most $1$. And $F(H)$ is a proper semigroup-ideal of $B(H)$ as the identity $I \not \in F(H)$.\\

 Regarding generators, the semigroup $B(H)$ is singly generated by the identity operator but $F(H)$ is not singly generated. Indeed, if it were generated by an operator $T$ of rank $k$, then $F(H) = \mathcal S(T) = \{T^n \mid n >0\}$ which is not possible because $\mathcal S(T)$ contains only elements of rank at most $k$.   
\eR

Another important distinction  between algebraic $B(H)$-ideals and multiplicative $B(H)$-semigroups and their multiplicative ideals that will be crucial in understanding the proofs of our main Theorems \ref{TR1} and \ref{TR2} is Remark \ref{scalars}: that semigroups may not be closed under scalar multiplication, e.g.,  for any  projection $P$, one has $P \in \SC(P) = \{P^n, n \ge 1\} = \{P\}$, but $2P \not \in \SC(P)$.
This substantially complicates the proof of the characterization of selfadjoint-ideal semigroups $\SC(T, T^*)$ generated by rank one $T$.

%\bR{R-SI}
%In summary, among all algebraic $B(H)$-ideals, $F(H)$ and $B(H)$ are the only semigroups that are SI but not simple. Indeed, the principal ideal generated by a rank one orthogonal projection is a proper semigroup-ideal in the semigroup $F(H)$ because every operator in that ideal is of rank at most $1$. And $F(H)$ is a proper semigroup-ideal of $B(H)$ as the identity $I \not \in F(H)$. 
%Moreover, the semigroup $B(H)$ is singly generated by the identity operator. But $F(H)$ is not singly generated because if it is generated by an operator $T$ of rank $k$, then $F(H) = \mathcal S(T) = \{T^n \mid n >0\}$ which is not possible because $\mathcal S(T)$ contains only elements of rank at most $k$.   
%\eR

\vspace{.2cm}

\textbf{SI unitary invariance versus similarity invariance}

The SI property for semigroups as well as simplicity are unitary invariance as proved in the following theorem. But as anticipated, SI is not similarity invariant (Remark \ref{example} which depends on Theorem \ref{TR2}).
\bT{T1}
Let $H_1$ and $H_2$ be Hilbert spaces.
For every semigroup $\SC \subset B(H_1)$ and  unitary operator $U \in B(H_1, H_2)$,
$\SC$ is an SI semigroup if and only if $U\SC U^{*}$ is an SI semigroup in $B(H_2)$. Additionally, $\SC$ is simple  if and only if $U\SC U^{*}$ is simple.
\eT
\bp{}
Suppose $\SC$ is an SI semigroup. 
It is clear that $U\SC U^{*}$ is also a semigroup.
If $J$ is a nonzero ideal of $U\SC U^{*}$ and $X \in J$, then we need to show that $X^{*} \in J$. 
For some $Y \in \SC,$ $X = UYU^{*}$.  
As $\SC$ is an SI semigroup, the principal ideal generated by $Y$ in $\SC$ is selfadjoint. Therefore $Y^{*} = Y_1 Y Y_2$ for some $Y_1, Y_2 \in \SC \cup \{I\}$.
Applying $U$ and $U^{*}$ to the latter equation yields $X^{*} = UY^{*}U^{*} =  UY_1 Y Y_2U^{*} = UY_1U^{*}UYU^{*}UY_2U^{*}$. 
Setting $X_1 := UY_1U^{*} $ and  $X_2:= UY_2U^{*}$, one has $X_1, X_2 \in U\SC U^{*}$, and since $J$ is an ideal of $U\SC U^{*}$, $X^{*} = X_1 X X_2 \in J$, showing $J$ is selfadjoint.
Therefore $U\SC U^{*}$ is an SI semigroup. The proof of the converse is symmetric and so clear.

The second part of the theorem follows directly from the fact that $J$ is an ideal in $\SC$ if and only if $UJU^{*}$ is an ideal in  $U\SC U^{*}$. 
\ep

The following Remark \ref{R3}(i)-(ii) and (iv)-(v) provide examples of simple and nonsimple SI semigroups and Remark \ref{R3}(iii) provides a necessary norm condition for a semigroup in $M_n(\mathbb C)$ or $B(H)$ to be an SI semigroup. 

\bR{R3}

(i) (Example of a nonsimple SI semigroup: $\SC(T)$ generated by certain selfadjoints) Consider a nonzero selfadjoint  operator $T$ with the property that $T^{n} \not \in \{I, T\}$ for all $n \geq 2$. (E.g., $T=\text{diag} \left<1/j\right>^{k}_{j=1}, \text{for any } 2 \leq k\leq \infty$.) The selfadjoint semigroup $\mathcal S(T)$ generated by $T$ consists of precisely all positive powers of $T$ and so contains only selfadjoint operators. This semigroup is not simple. In fact, it properly contains a family of  principal ideals: $(T^{m})_{\SC(T)}$ for each $m \geq 2$. Indeed, if $(T^{m})_{\SC(T)} = \SC(T)$, then $T^{k} = T$ for some $k \geq 2$ which violates the condition $T^{n} \not \in \{I, T\}$ for each $n \geq 2$.  Also $\SC(T)$ is clearly an SI semigroup because $\SC(T)$ contains only selfadjoint operators so all its nonzero ideals are selfadjoint.

(ii) (Examples of simple and nonsimple SI semigroups: $\SC(T)$ generated by certain selfadjoints)

(a) An obvious example for simple is to choose $T$ a projection operator, then $\SC(T) = \{T\}$ is simple. A less trivial class is to choose $T$ a nonzero selfadjoint operator with $T^n = T$ for $n >2$. Then $\SC(T) = \{T, T^2, \cdots, T^{n-1}\}$ and note that this selfadjoint semigroup does not contain any proper nonzero ideal because for $J$ a nonzero ideal and $0 \neq T^{j} \in J$, $T = T^n = T^{n-j}T^{j} \in J$ for $1 < j \leq n-1$ implying $\SC(T) \subset J$ hence $J = \SC(T)$, so $\SC(T)$ is simple.

(b) For $0 <a <1$ and \begin{equation*}
T = \begin{pmatrix}
a&0\\
0&0
\end{pmatrix},
\end{equation*}
the semigroup $\SC(T) = \{T^n\mid n \geq 1\}$ is not simple. Indeed, $(T^2)_{\SC(T)} \neq \SC(T)$ because $T \not \in (T^2)_{\SC(T)} = \{T^n \mid n\geq 3\}$.

(iii) (Necessary norm condition for SI property when containing a nonselfadjoint operator)\\ If a semigroup $\SC$ contains a nonselfadjoint operator $T$, then a necessary condition for $\SC$ to be an SI semigroup is that $\SC$ must contain at least one element of norm greater than or equal to $1$. Indeed, consider the principal ideal $(T)_{\SC}$ of $\SC$. For this ideal to be selfadjoint, $T^{*} = XTY$ for $X, Y \in \SC \cup {\{I\}}$ and $X$ or $Y \in \SC \setminus \{I\}$. So  $0<||T||= ||T^{*}|| = ||XTY|| \leq ||X||||T||||Y||$ and so either $||X|| \geq 1$ or $||Y|| \geq 1$.

(iv) Selfadjoint semigroups of unitaries (which are in fact groups) are simple.  Indeed, if $T\in J$ an ideal of selfadjoint semigroup $\SC$ containing only unitaries, then $T^* \in \SC$ and $T^*T = I \in J$ implying $J = \SC$.

(v) For $V$ an isometry (or co-isometry), $\SC(V, V^*)$ is simple. Indeed, as $V^*V = I$ (or $VV^* = I$), viewed abstractly $\SC(V, V^*)$  is a bicyclic semigroup, and the list for $\SC(T, T^*)$ in \text{Proposition} ~\ref{P1} reduces to  $$ \{I, V^n, {V^*}^m, V^n{V^*}^m \mid n,m >0\}.$$ 
Using the identity $V^*V = I$, one has ${V^*}^nV^n = I$, ${V^*}^mV^m = I$, and ${V^*}^nV^n{V^*}^mV^m= I$ for $m, n > 0$. This implies that the identity operator belongs to all the principal ideals $(V^n)_{\SC(V, V^*)}, ({V^*}^m)_{\SC(V, V^*)}$, and $(V^n{V^*}^m)_{\SC(V, V^*)}$. So for $m,n >0$,
$$(V^n)_{\SC(V, V^*)}= ({V^*}^m)_{\SC(V, V^*)} = (V^n{V^*}^m)_{\SC(V, V^*)}=\SC(V, V^*).$$
Hence all nonzero principal ideals coincide with $\SC(V, V^*)$. Therefore, any nonzero ideal of $\SC(V, V^*)$ coincides with $\SC(V, V^*)$ since every nonzero ideal contains a nonzero principal ideal.
\eR
The necessary norm condition of Remark \ref{R3}(iii) above enables us to construct obvious examples of non-SI semigroups as given below.

\begin{example}\label{EE}
The selfadjoint semigroup $ \mathcal S(T, T^{*})$ with $T$  a nonselfadjoint operator with $||T|| < 1$ is not an SI semigroup because every element of $\mathcal S(T, T^*)$ has norm strictly less than $1$, hence violating the necessary norm condition Remark \ref{R3}(iii).  
\end{example}
 
\newpage 
\section{SI semigroups $S(T, T^{*})$ for normal operators and partial isometries}

Semigroups in $B(H)$ with the SI property being a new class of operators to study,  at this point one is clueless about the composition of this special class of semigroups. For instance, whether SI semigroups in $B(H)$ are sparse or abundant, and what classes of operators have their $\SC(T, T^*)$ an SI semigroup? A natural place to start the investigation is to consider selfadjoint semigroups generated by sets. In this paper we focus on the study of singly generated selfadjoint semigroups, semigroups $\mathcal S(T, T^{*})$ generated by $\{T, T^*\}$, because as we shall see, this is a rich subject.

Naturally the semigroup $\SC(T, T^*)$ generated by a selfadjoint operator $T$ is always an SI semigroup (Remark \ref{sa}(ii)), but characterizing semigroups generated by nonselfadjoint operators with the SI property is much more complicated. We shall investigate the SI property of $\SC(T, T^*)$ first for normal nonselfadjoint operators, perhaps the most natural nonselfadjoint operators to consider. We determine those normal nonselfadjoint  operators $T$ for which $\SC(T, T^*)$ is an SI semigroup.  And for these $T$,  $\SC(T, T^*)$ turns out to be vacuously SI, i.e., simple.

\bT{N2}
For $T$ a normal nonselfadjoint operator, the following are equivalent.
\begin{enumerate}[label=(\roman*)]
\item $\SC(T, T^*)$ is an SI semigroup.
\item $T$ is unitarily equivalent to $U \oplus 0$ (or $U$ when $\text{ker} T= \{0\}$) with $U$ a unitary operator. 
\item $\SC(T, T^*)$ is a simple semigroup.
\end{enumerate}
\eT

\bp
By the spectral theorem for normal operators \cite[Theorem 1.6]{RR1}, $T$ is unitarily equivalent to a multiplication operator $M_{\phi}$ on $L^{2}(\textbf{X}, \mu)$ with $\phi$ an essentially bounded measurable function on a finite measure space $(\textbf{X}, \mu)$. Since $\SC(T, T^*)$ is an SI semigroup,
it follows from Theorem \ref{T1} that $\mathcal S(M_{\phi}, M_{\bar{\phi}})$ is also an SI  semigroup. Since $T$ is a nonselfadjoint operator, $M_{\bar{\phi}} \neq M_{\phi}$.

The elements of $\mathcal S(M_{\phi}, M_{\bar{\phi}})$ consists of all words in $M_{\phi}$ and $M_{\bar{\phi}}$ (recall paragraph following Definition \ref{D6}), but $M_{\phi}$, $M_{\bar{\phi}}$ commuting  together with Proposition \ref{P1} provides a simpler description: 
 $$\mathcal S(M_{\phi}, M_{\bar{\phi}}) = \{M_{{\phi}^n}, M_{\bar{\phi}^{m}}, M_{\phi^{n}\bar{\phi}^{m}} \mid m, n \geq 1\}.$$ 

(i)$\Rightarrow$(ii): Suppose $\SC(T, T^*)$ is an SI semigroup. So the principal ideal $(M_{\phi})_{\mathcal S(M_{\phi}, M_{\bar{\phi}})}$ of $\mathcal S(M_{\phi}, M_{\bar{\phi}})$ is selfadjoint. This implies $M_{\bar{\phi}} = XM_{\phi}Y$  where $X, Y \in \SC(M_{\phi}, M_{\bar{\phi}}) \cup {\{I\}}$ and $X$ or $Y \in \SC(M_{\phi}, M_{\bar{\phi}}) \setminus \{I\}$ where $I$ is the identity operator (apply Lemma \ref{L2} (i)$ \Rightarrow$(iii) and since $M_{\bar{\phi}} \neq M_{\phi}$, $X$ or $Y\neq I$).
 By display above, since $M_\phi \neq M_{\bar{\phi}}$, one must have 
 $M_{\bar{\phi}} = M_{\phi^{n}}$ for $n > 1$ or if not, then $M_{\bar{\phi}}= M_{\phi^{n}\bar{\phi}^{m}}$ for some $n, m \geq 1$ since $M_{\bar{\phi}} = XM_{\phi}Y$ and $XM_{\phi}Y$ contains in the semigroup product one $M_{\phi}$ along with at least one other $M_{\phi}$ or $M_{\bar{\phi}}$. (Note: besides $M_{\bar{\phi}} = XM_{\phi}Y$, $M_{\bar{\phi}}$ could simultaneously could also be equal to other forms from the displayed list above.) Hence either $\bar{\phi}(x) = \phi^{n}(x)$ a.e. or $\bar{\phi}(x) = \phi^{n}(x)\bar{\phi}^{m}(x)$ a.e., and so $|\phi| = |\phi|^{k}$ a.e. for some $k > 1.$
For $E := \{x \in \textbf{X} \mid |\phi(x)| = |\phi(x)|^{k}\}$,  $\mu(E^{c}) = 0.$ 
Therefore $|\phi(x)|$ is either $0$ or $1$ on $E$. Let $F:= \{x \in E \mid |\phi(x)| = 1\}$.
Since $T$ is nonselfadjoint, $T \neq 0$ and hence $\mu(F) >0$.
Splitting the Hilbert space $L^2(\textbf{X}, \mu) = L^2(F, \mu|_{F})\oplus L^2(\textbf{X}\setminus F, \mu|_{\textbf{X} \setminus F})$, then $M_{\phi}$ splits as $M_{\phi} =  M_{\phi \chi_{F}} \oplus M_{\phi\chi_{\textbf{X}\setminus F}} =  M_{\phi \chi_{F}} \oplus 0$  and $M_{\phi \chi_{F}}$ is a unitary operator on $L^2(F, \mu|_{F})$ as $|\phi| = 1$ on $F$.  
Therefore $T$ is unitarily equivalent to $U \oplus 0$ where $U$ is a unitary operator. That is (i) $\Rightarrow$ (ii). This completes the proof of Theorem \ref{N2}. 

Additionally there is more structure information for $T$.
Indeed, the above phrase when $\SC(T, T^*)$ is an SI semigroup ``either $\bar{\phi}(x) = \phi^{n}(x)$ a.e. or $\bar{\phi}(x) = \phi^{n}(x)\bar{\phi}^{m}(x)$ a.e." has two interesting further implications which however will not be used elsewhere in this paper. Firstly, as just shown, this pair of conditions implies $|\phi | = \chi_F$ a.e.. But $|\phi | = \chi_F$ a.e. is equivalent to $\bar{\phi} = \phi \bar{\phi}^2$. So the operator theoretic conditions led us to this pair of conditions which are equivalent to the single condition $|\phi | = \chi_F$ a.e..  However, the pair of conditions has a noncommutative analog essential for our main theorem on SI characterization of $\SC(T, T^*)$ for non-normal rank one operators (equivalently, non-nilpotent rank one operators). (See Theorem \ref{TR2} (ii) trace-norm condition.)

Secondly, if in addition to the condition $|\phi | = \chi_F$ a.e. either $\bar{\phi}(x) = \phi^{n}(x)$ a.e. or $\bar{\phi}(x) = \phi^{n}(x)\bar{\phi}^{m}(x)$ a.e. for $n,m \geq 1$ and $n+1 \neq m$, then $\phi(F) \subset \mathbb{S}^1$ (the unit circle),  and furthermore, on $F$ in both cases it follows that $\phi^{k}(x) \equiv 1$ for some $k \geq 1$. This follows after multiplying each equation by $\phi$. Then defining $E_n:= \{x \in X | \,\phi(x) = \exp(\frac{2\pi i n}{k})\}$ for $1 \leq n \leq k$, $X$ partitions into $X = \cup_{n}E_n \cup \{X \setminus \cup_{n}E_n\}$. So one has the finite decomposition $$L^2(X, \mu) = \oplus_nL^2(E_n, \mu|_{E_n}) \oplus L^2(X\setminus \cup_{n}E_n) , \mu|_{X\setminus \cup_{n}E_n})$$ on which $M_{\phi} = \sum^{\oplus}\{M_{\alpha_nI} \mid \mu(E_n) >0\} \oplus 0$ where $\alpha_ n := \exp(\frac{2\pi i n}{k})$ for $1 \leq n \leq k$. Therefore $T$ is unitarily equivalent to $\sum^{\oplus}\{M_{\alpha_nI} \mid \mu(E_n) >0\} \oplus 0$ where $\{\alpha_n\}_{n = 1}^{k}$ is the set of all $k^{th}$ roots of unity and where the $0$ summand may not appear.\\

\noindent (ii)$\Rightarrow$(iii): By Theorem \ref{T1}, it suffices to prove for $T = U \oplus 0$ with $U$ a unitary operator that $\SC(T, T^*)$ is a simple semigroup. 
Since $UU^{*} = U^{*}U = I$, $\SC(T, T^{*}) = \{U^{n} \oplus 0, U{^{*}}^{n} \oplus 0, I \oplus 0 \mid n \geq 1\}$ by Proposition \ref{P1}, which by inspection is a group. For every $X \in J$ an ideal of $\SC(T, T^{*})$, $X \in \SC(T, T^*)$ and hence for some $n \geq 1$, either $X = U^{n} \oplus 0$, $X = U{^{*}}^{n} \oplus 0$  or $X = I \oplus 0$. Since $J$ is an ideal, for all three cases, $X^{*}X = I \oplus 0 \in J$, so $J$ contains $I \oplus 0$, the identity of the semigroup $\SC(T, T^{*})$, so $\SC(T, T^{*}) \subset J$, i.e., $J = \SC(T, T^{*})$ and is in fact a group. Therefore $\SC(T, T^{*})$ contains no nonzero proper ideal, i.e., simple. 

\noindent (iii)$\Rightarrow$(i): As simple semigroups have no proper nonzero ideals, so $\SC(T, T^*) $ is vacuously an SI semigroup.
\ep

The discussion in the second paragraph of this section combined with Theorem \ref{N2} provides a characterization of singly generated SI semigroups $\SC(T, T^*)$ for $T$ normal.
\bC{T2}
(Structure theorem for singly generated SI semigroups $\SC(T, T^*)$ for $T$ normal.)
\begin{enumerate}[label=(\roman*)]
\item If $T$ is selfadjoint, then $\SC(T, T^*) = \{T^n, n >0\}$ is an SI semigroup (not necessarily simple, see Remark \ref{R3}(i)-(ii)). 
\item If $T$ is a normal nonselfadjoint operator, then $\SC(T, T^*)$ is simple.
\end{enumerate} 
\eC

 \bR{R4}
The operator $T = U \oplus\, 0$ ($U$ unitary) that appeared in Theorem \ref{N2}(ii) belongs to the class of power partial isometries (\cite[Theorem]{HW}, also recall seventh paragraph in the Introduction section), so Theorem \ref{N2} links SI semigroups generated by a normal nonselfadjoint operator to semigroups generated by a power partial isometry and simple selfadjoint semigroups. 
 The Popov-Radjavi fact that an operator $T$ is a power partial isometry if and only if the selfadjoint semigroup $\SC(T, T^*)$ generated by $T$ consists of only partial isometries \cite[Proposition 2.2]{RP} combined with our Corollary \ref{CPP} herein, together imply: 
 \begin{center}
 \textit{Singly generated selfadjoint semigroups $\SC(T, T^*)$ generated by a power partial isometry $T$ are SI.} 
 \end{center}
But $\SC(T, T^*)$ generated by a power partial isometry may not be simple (Example  \ref{ENS}). 
\vspace{.3cm}

Characterizing which of these are simple is an interesting problem. First recall that every power partial isometry is a direct sum $T = \mathcal U \oplus \mathcal V \oplus \mathcal W \oplus \mathcal S$ whose summands are direct sums of unitary operators, pure isometries, pure co-isometries, and truncated shifts with the understanding that not all summands must necessarily be present \cite[Theorem]{HW}. This suggests different types of partial isometries that might be used for this characterization. For instance, we have already proved that the selfadjoint semigroup $\SC(T, T^*)$ is simple if $T$ is a unitary, a pure isometry (or a pure co-isometry), a rank one non-normal power partial isometry (Theorem \ref{N2},  Remark \ref{R3}(v), and Corollary \ref{CPI} respectively). (Recall that a pure isometry is just a unilateral shift of arbitrary multiplicity, a pure co-isometry is the adjoint of a pure isometry, and a truncated shift is a Jordan block with zero diagonal, i.e., a truncated shift $S$ of index $n\geq 1$ is defined by taking an $n$-fold direct sum $H \oplus \cdots \oplus H$ and defining $S =0$ when $n=1$ and $S(f_1, f_2, \cdots, f_n) := (0, f_1, \cdots, f_{n-1})$ when $n>1$.)  \eR

\bE{ENS} For $T$ a power partial isometry, the necessarily SI semigroup $\SC(T, T^*)$ is \underline{not simple} if $T$ is of the form $T =  U \oplus V \oplus W \oplus S$ with $\text{index}(S) >1$, where $U^*U = UU^* = V^*V = WW^* =I$.
 Indeed, it suffices to construct a proper nonzero ideal of $\SC(T, T^*)$ to prove non-simplicity of $\SC(T, T^*)$.  Since for $n = \text{index(S)} >1$, one has $S^n \equiv 0$ but $S \not \equiv 0$. Then the operator $A= U^n \oplus V^n \oplus W^n \oplus 0 \in \SC(T, T^*)$. By applying Lemma \ref{L1} to $\mathcal A = \{A\}$, 
 $$(A)_{\SC(T, T^*)} = \{A\} \cup  \SC(T, T^*)\{A\} \cup \{A\}\SC(T, T^*) \cup \SC(T, T^*) \{A\} \SC(T, T^*).$$
 Observe that $0$ appears in the last direct summand of every operator in the principal ideal $(A)_{\SC(T, T^*)}$, because $A$ and all elements of $\SC(T, T^*)$ are simultaneously block diagonal ending in block $0$. So $T = U \oplus V \oplus W \oplus S  \not \in  (A)_{\SC(T, T^*)}$ implying $(A)_{\SC(T, T^*)} \subsetneq \SC(T, T^*)$, i.e., $\SC(T, T^*)$ is not simple.
 \eE
 An example of a partial isometry $T$ for which the semigroup $\SC(T, T^{*})$ is not an SI semigroup is: 

\bE{E}{\textit{(A non-SI semigroup $\SC(T, T^*)$ in $M_{2}(\mathbb C)$ generated by a partial isometry)}}\\
For $0 <|a|<1,\, 0 <|b|<1$ and $|a|^2 +|b|^2 = 1$ with $a \neq \bar{a}$, the $2 \times 2$ matrix
 \begin{equation*}
T = \begin{pmatrix}
a&b\\
0&0
\end{pmatrix}.
\end{equation*}
 is a partial isometry (since $T = TT^*T$, a partial isometry characterization \cite[Corollary 3 of Problem 127]{H}), and  $T$ is not a power partial isometry ($T^2$ is not a partial isometry because $||T^2|| = |a| <1$ whereas nonzero partial isometries have norm $1$). But $\SC(T, T^*)$ is not an SI semigroup. Indeed, this is a special case covered by Theorem \ref{TR2} when $T$ is a rank one operator in $M_{2}(\mathbb C)$ with nonzero trace. That is, Theorem \ref{TR2} condition (ii) fails because $tr T$ is not real and the trace-norm condition fails because for this $T$, $||T|| =1$ and so the trace-norm condition implies ${tr^m T}\overline{tr ^nT} = 1$ for some $m,n \geq 1$, hence $|tr T| =1$ against our assumption $|tr T| <1$. 
 \eE
 Summarizing Theorem \ref{N2}, Remark \ref{R4}, and Examples \ref{ENS}-\ref{E}, for normal $T$, $\SC(T, T^*)$ is SI if and only if $T = U \oplus 0$ (a power partial isometry) if and only if $\SC(T, T^*)$ is simple; generally for $T$ a partial isometry $\SC(T, T^*)$ may not be SI; for $T$ a power partial isometry of the form $T = U \oplus V \oplus W \oplus S$, $\SC(T, T^*)$ is SI but when the truncated shift $S \neq 0$, $\SC(T, T^*)$ is not simple; and when $T = U \oplus 0$, $\SC(T, T^*)$ is simple.\\

Examples \ref{ENS}-\ref{E} lead naturally to the possibility of characterizing nonsimple SI semigroups $\SC(T, T^*)$ and simple semigroups  $\SC(T, T^*) $ generated by partial isometries $T$ in Question \ref{QPI} below. 

\begin{question}\label{QPI}
Characterize which partial isometries $T$ have SI semigroup $\SC(T, T^*)$, and among those determine which are simple.
\end{question}

  We first have a complete answer to the SI and simplicity of $\SC(T, T^*)$ for isometries and co-isometries.
 Indeed, we proved in Remark \ref{R3}(v) that for isometries and co-isometries $T$, $\SC(T, T^*)$ is simple.\\

For rank one non-normal partial isometries, we have Corollary \ref{CPI}:
\begin{enumerate}[label=(\roman*)]
\item if $T$ is a power partial isometry, then $\SC(T, T^*)$ is simple.
\item if $tr T \in \mathbb R\setminus \{0\}$, then $\SC(T, T^*)$ is an SI semigroup.
\end{enumerate}
But for higher ranks, $T$ a partial isometry not a power partial isometry we do not have a characterization for $\SC(T, T^*)$ simple nor for $T$ a power partial isometry of the form $T = \mathcal U \oplus \mathcal V \oplus \mathcal W \oplus 0$ where $\mathcal U, \mathcal V$ and $\mathcal W$ are classes as discussed above and where not all summands must necessarily be present. 
\vspace{.2cm}

\section{Characterization of SI and simple semigroups $S(T, T^{*})$\\ for rank one operators in $B(H)$}

In this section, we completely characterize the SI semigroup $\SC(T, T^*)$ for $T$ rank one by dividing them into the two cases, $trT= 0$ and $trT \neq 0$, each with its own characterization. At the end of the section, we summarize the characterizations of SI semigroups $\SC(T, T^*)$ for $T$ normal and $T$ rank one leading to the further classification of these SI semigroups into simple and nonsimple semigroups.\\

As this section concerns rank one operators, we first recall some tensor facts:\\ For $f, g \in H$, $(f \otimes g)(h):= (h, g)f$ represents all possible rank one operators in $B(H)$. \\Some standard relations for vector tensor products are:
\begin{enumerate}[label=(\roman*)]
\item $(f \otimes g)(f' \otimes g') = (f', g) (f \otimes g')$ as $$(f \otimes g)(f' \otimes g')(h) = (f \otimes g)(h, g')f' = (h, g')(f', g)f = (f', g) (f \otimes g')(h).$$

\item $(f \otimes g)^* = g \otimes f$ as $(h, (f \otimes g)^*k) = ((f \otimes g)h, k) = (h, g)(f, k)$ and $$(h, (g\otimes f)k = \overline{(k,f)}(h,g) = (h, g)(f, k).$$
\item $||f \otimes g|| = ||f||||g||$ as $||f \otimes g|| = sup_{||h|| =1} \{||(f \otimes g)h||\} = sup_{||h|| =1} \{|(h,g)| ||f||\} = ||g||||f||$. 
\end{enumerate}
\vspace{.2cm}

\noindent And when $||f|| =1$, $f \otimes f$ is a rank one projection. \\

A useful matrix fact for rank one operators $T$ is their unitary equivalence:
 \begin{equation}\label{Mm}
T \cong \begin{pmatrix}
a&b\\
0&0
\end{pmatrix} \oplus 0
\text{, \quad for some ~} a, b \in \mathbb C, \text{where $\oplus~ 0$ may be absent.} \end{equation} 
\noindent This observation will occasionally  help simplify some computations, but in general we will employ abstraction in part in the hope of applying this section to higher ranks. For instance, see trace-norm equalities (Theorem \ref{TR2}), tensors as described above, and relations (e.g. Proposition \ref{PT3}). 
\bP{PT1}
For $T$ a rank one operator in $B(H)$, let $Ran T = span \{f\}$ normalized by $||f|| = 1$. Then $T = f \otimes T^*f$, $T^*f \neq 0$, $T^n = (Tf,f)^{n-1}T$ for $n \ge 2$, $tr T = (Tf,f)$ and $||T|| = ||T^*f||$. Consequently, $T^n = (tr T)^{n-1}T$.  
\eP
\bp
For each $h \in H$, $Th = af$ for some $a \in \mathbb C$. So $(Th,, f) = (af, f) = a$, and hence 
$$Th = (Th, f)f = (h, T^*f)f = (f \otimes T^*f)h$$ for all $h \in H$, i.e., $T = f \otimes T^*f \neq 0$, and hence $T^*f \neq 0$. Using relation $(f\otimes g)(h \otimes k) = (h, g)(f \otimes k)$, one obtains $T^n = (Tf, f )^{n-1}T$ for $n \ge 2$.  For $\{g_n\}$ an orthonormal basis for $Ran^{\perp} T$, $\{f\} \cup \{g_n\}$ forms an orthonormal basis for $H$. Therefore $tr T = (Tf, f ) + \sum_n(Tg_n, g_n) = (Tf, f )$ because all $g_n \perp Ran T$. And using Property (iii) of rank one tensor products, $||T|| = ||f \otimes T^*f|| = ||f||||T^*f|| = ||T^*f||$. 
\ep

For the case: non-normal rank one with $tr T \neq 0$, we need the following proposition used later in the proof of Theorem \ref{TR2}, where the condition $Ran T^* \not \subset Ran T$ is essential.

\bP{PN}
A rank one operator $T$ is normal if and only if $Ran T^* \subset Ran T$.
\eP
\bp
Here we use the matrix form for $T$ (Equation (\ref{Mm}) mentioned above).
Since $TT^* = T^*T$,
\begin{equation*}\label{M1}
\begin{pmatrix}
|a|^2+|b|^2 &0\\
0&0
\end{pmatrix} \oplus 0 = \begin{pmatrix}
a&b\\
0&0
\end{pmatrix}  \begin{pmatrix}
\bar{a}&0\\
\bar{b}&0
\end{pmatrix} \oplus 0  = \begin{pmatrix}
\bar{a}&0\\
\bar{b}&0
\end{pmatrix} \begin{pmatrix}
a&b\\
0&0
\end{pmatrix} \oplus 0= \begin{pmatrix}
|a|^2 & \bar{a}b\\
a\bar{b} & |b|^2
\end{pmatrix} \oplus 0
\end{equation*}
So $T$ is normal if and only if $b=0$ if and only if $Ran T^* = Ran T$. Moreover, as $T$ is rank one, so $b=0$ and $a \neq 0$ is equivalent to the weaker condition $Ran T^* \subset Ran T$.
\ep
\bR{RN}
This proof shows the stronger result that normality of a finite-rank operator implies $RanT^*= Ran T$, but the converse fails easily at higher ranks. For instance, any invertible non-normal operator has $RanT^*= Ran T$.\eR

\bP{PT3}
For rank one operators $T \in B(H)$ and $n \geq 1$,
\begin{enumerate}[label=(\roman*)]
\item $(T^*T)^n= ||T||^{2n-2}T^*T$ \text{and} $(T^*T)^{2n} = ||T||^{2n}(T^*T)^n$
\item $(T^*T)^nT^* = ||T||^{2n}T^*$
\item $(TT^*)^n = ||T||^{2n-2}TT^*$ \text{and} $(TT^*)^{2n} = ||T||^{2n}(TT^*)^n$
\item $(TT^*)^n T= ||T||^{2n}T$
\item $T^n = (tr T)^{n-1}T$ when $n > 1$.
\end{enumerate}
\eP
\bp
Substituting $f \otimes T^*f$ for $T$ and $T^*f \otimes f$ for $T^*$ with $||T|| = ||T^*f||$, $||f|| =1$, and using $(f\otimes g)(h \otimes k) = (h, g)(f \otimes k)$ repeatedly, one obtains (ii) via
\begin{align*}
(T^*T)^nT^* &= T^*(TT^*)^n 
= (T^*f\otimes f)[(f\otimes T^*f)(T^*f\otimes f)]^n 
= (T^*f\otimes f)(||T^*f||^2(f\otimes f))^n \\&=  ||T^*f||^{2n}(T^*f\otimes f)(f \otimes f) = ||T^*f||^{2n}T^* = ||T||^{2n}T^*.
\end{align*}

\noindent Using (ii) yields (i) (first part) via $(T^*T)^n = (T^*T)^{n-1}(T^*T) = [(T^*T)^{n-1}T^*]T =  ||T||^{2n-2}T^*T$ when $n \geq 1$.

One obtains (iv) directly via
\begin{align*}
(TT^*)^nT &= [(f \otimes T^*f)(T^*f \otimes f)]^nT = [||T^*f||^2f \otimes f]^nT\\
&=  ||T||^{2n}(f \otimes f)(f \otimes T^*f ) =||T||^{2n} f \otimes T^*f = ||T||^{2n}T
\end{align*}

Similarly, one obtains (iii) (first part)  via $(TT^*)^n = (TT^*)^{n-1}TT^* = [(TT^*)^{n-1}T]T^* = ||T||^{2n-2}TT^*$.\\

Hence we have (i)(second part) from (i) (first part) via
\begin{align*}
(T^*T)^{2n} &= ||T||^{4n-2}T^*T = ||T||^{2n}||T||^{2n-2}T^*T = ||T||^{2n}(T^*T)^n
\end{align*}
and (iii) (second part) from (iii) (first part) via
\begin{align*}
(TT^*)^{2n} &=||T||^{4n-2}TT^*= ||T||^{2n}||T||^{2n-2}TT^*=  ||T||^{2n}(TT^*)^n.
\end{align*}

Also (v), $T^n = (tr T)^{n-1}T$, follows directly from Proposition \ref{PT1}.
\ep

Moreover, one can identify rank one partial isometries by their norm via
\bP{PT2}
For rank one $T \in B(H)$, $T$ is a partial isometry if and only if $||T|| = 1.$
\eP

\bp
This is an instance where a matrix proof is simpler than an abstract proof like the one above, and gives a bit more information used later in Proposition \ref{Ptr}.

Recall the condition that characterizes partial isometries: $T = TT^*T$ \cite[Corollary 3 of Problem 127]{H}. Using the matrix form for $T$ obtains  \begin{equation*}
T = \begin{pmatrix}
a&b\\
0&0
\end{pmatrix} \oplus 0 = TT^*T = \begin{pmatrix}
a(|a|^2 + |b|^2) & b(|a|^2 + |b|^2)\\
0&0
\end{pmatrix} \oplus 0.
\end{equation*}
Therefore, $T$ is a partial isometry if and only if $|a|^2 + |b|^2 =1$, i.e., $||T|| = \sqrt{|a|^2 + |b|^2} = 1$.

\ep

%\bp
%By Proposition \ref{PT1}, $T = f \otimes T^*f$  for some $f$ with $||f|| =1$.
%As stated earlier, a condition to characterizing  partial isometries $T$ is given by the identity $T = TT^*T$ \cite[Corollary 3 of Problem 127]{H}.
%Substituting $0 \neq T = f \otimes T^*f$, one obtains $$f \otimes T^*f = (f \otimes T^*f)(T^*f \otimes f)(f \otimes T^*f) = ||T^*f||^2 (f \otimes T^*f) = ||T||^2(f \otimes T^*f),$$ and so $f \otimes T^*f = ||T||^2(f \otimes T^*f)$ which is equivalent to $||T|| = 1$ because $T = f \otimes T^*f \neq 0$.

\bP{PI2}
For a rank one partial isometry $T \in B(H)$, 
\begin{center}
$T$ is a power partial isometry if and only if either $tr T = 0$ or $|tr T| = 1$. 
\end{center}
\eP
\bp
Suppose $T$ is a power partial isometry and $tr T \neq 0$. Then for each $n\geq1$, $T^n = {T}^n{T^*}^n{T}^n$. And for $n > 1$, $T^n = (tr T)^{n-1}T$ (Proposition \ref{PT3}(v)). So $$(tr T)^{n-1}T = {T}^n{T^*}^n{T}^n=|tr T|^{(2n-2)}(tr T)^{n-1} TT^*T =  |tr T|^{(2n-2)}(tr T)^{n-1}T.$$
Since $T \neq 0$ (as $T$ is rank one) and $tr T \neq 0$, so for $n>1$, $|tr T| =1$. Therefore $|tr T| = 1$.

Conversely, suppose $|tr T| =1$. As $T^n = (tr T)^{n-1}T$ for $n \geq 1$, so  $(tr T)^{n-1}(\overline{tr T})^{n-1} =1$ and
 $$(tr T)^{n-1}T = T^n = (tr T)^{n-1}(\overline{tr T})^{n-1}(tr T)^{n-1} T = (tr T)^{n-1}(\overline{tr T})^{n-1}(tr T)^{n-1}TT^*T = T^n{T^*}^nT^n.$$
Therefore $T^n = {T}^n{T^*}^n{T}^n$, i.e., $T$ is a power partial isometry.

Finally if $tr T = 0$, then for $n >1$,  
$$T^n = (tr T)^{n-1}T = 0 = T^n{T^*}^nT^n.$$
For $n=1$, $T$ is a partial isometry as it is a part of the hypothesis. 
Therefore, $T$ is a power partial isometry. 
\ep

\noindent \textit{Observation}. The rank one power partial isometry conditions in Proposition \ref{PI2}, $tr T = 0$ or $|tr T| = 1$, are equivalent to $T$ being nilpotent or a modulus one multiple of a projection. Indeed, observe from Equation (\ref{Mm}) that $a = tr T$ so $tr T =0$ implies that $T$ is nilpotent. And when $|tr T| =1$ then $b=0$ as $||T|| =1$ by Proposition \ref{PT2}. So $|a| = |tr T| = 1$ and so $T$ is a modulus one multiple of a projection. 

Another connection between power partial isometries and nilpotence comes from the Halmos-Wallen characterization of power partial isometries: $T \cong U \oplus V \oplus W \oplus S$ which clearly is nilpotent if and only if $U = V= W =0$ and the order of nilpotence is the index of $S$ (see Example \ref{ENS} and preceeding paragraph).\\

Working with rank one operators has the advantage that it simplifies somewhat the complicated expressions in $T$ and $T^*$ that occur in the description of the semigroup $\SC(T, T^*)$ listed in Proposition \ref{P1}. We describe increasingly simplified forms for $\SC(T, T^*)$ in Remark \ref{RS} below. 
 \newpage
\bR{RS}{(Simpler forms of $\SC(T, T^*)$ for $T$ a rank one operator)}\\
For $T$ a rank one operator (so nonzero), the Proposition \ref{P1} list for $\SC(T, T^*)$ reduces to simpler forms as follows.\\
Recall this list is:
\begin{equation}\label {E1}
 S(T, T^*) = \{T^n, {T^*}^n, \Pi_{j=1}^{k}T^{n_j}{T^*}^{m_j},  (\Pi_{j=1}^{k}T^{n_j}{T^*}^{m_j})T^{n_{k+1}}, \Pi_{j=1}^{k}{T^*}^{m_j}T^{n_j}, (\Pi_{j=1}^{k}{T^*}^{m_j}T^{n_j}){T^*}^{m_{k+1}}\} 
 \end{equation}
 where $n \ge 1,\,  k\ge1,\, n_j, m_j \ge 1\, \text{for}\, 1 \le j \le k, ~\text{and}~ n_{k+1}, m_{k+1} \geq 1$. The products $\Pi_{j= 1}^{k}$ denote of course ordered products. \\
 
Recall Proposition \ref{PT3}(v): $T^n = (tr T)^{n-1}T$ for $n >1$. For convenience of notation set $a := tr T = (Tf, f)$ so $T^n = a^{n-1}T$.\\
 
 The simpler forms are:

(i) If $a = 0$, then $T^2 = (tr T)T = 0$ and so $0 \in \SC(T, T^*)$ and List (\ref{E1}) for the semigroup $S(T, T^*)$ reduces to
\begin{equation} \label{E2}
\SC(T, T^*) = \{0, T, T^*, (TT^*)^n, (TT^*)^nT, (T^*T)^n, (T^*T)^nT^*\mid n\ge 1\}
\end{equation}
as the only possible nonzero operators occurring in $S(T, T^*)$ are products of terms alternating between $T$ and $T^*$ including singletons $T$ and $T^*$ alone. And using Proposition \ref{PT3}, the list reduces to the lower order terms:
\begin{equation}\label{E3}
\SC(T, T^*) = \{0, T, T^*, ||T||^{2n-2}TT^*, ||T||^{2n}T, ||T||^{2n-2}T^*T, ||T||^{2n}T^*\mid  n\ge 1\}.
\end{equation}
And from this list one can see why each of the $6$ types of possibly nonzero operators is in fact nonzero.\\

(ii) The case $a \neq 0$ reduction for List (\ref{E1}) is more complicated. If $a \neq 0$, then $T^n = a^{n-1}T$ for $n \ge 1$, using the convention that $a^0 = 1$. For example, the List (\ref{E1}) operator $\Pi_{j=1}^{k}T^{n_j}{T^*}^{m_j} \in S(T, T^*)$ simplifies 
to $a^{n_1-1}\bar{a}^{m_1-1}\cdots a^{n_k-1}\bar{a}^{m_k-1}(TT^*)^k$, so by rearranging the powers of $a$ and $\bar{a}$ one has the product $\Pi_{j=1}^{k}T^{n_j}{T^*}^{m_j} = a^p\bar{a}^q(TT^*)^k$ where $p,q \geq 0$ and $k\geq1$. Hence the containment for List (\ref{E1}): $$S(T, T^*) \subset \{a^nT, \bar{a}^nT^*, a^p\bar{a}^q(TT^*)^k, a^p\bar{a}^q(TT^*)^k T, a^p\bar{a}^q(T^*T)^k,  a^p\bar{a}^q(T^*T)^kT^* \mid n\ge0, p, q \ge 0, k\ge 1\}.$$
It is important to observe that the list on the right has scalars while List (\ref{E1}) on the left does not.

For the reverse inclusion note that, because $T^n = a^{n-1}T$ for $n \ge 2$, one has $T^{n+1}= a^nT \,\text{and\,}\\ {T^*}^{n+1} = \bar{a}^nT^*  \in \SC(T, T^*)$, $a^p\bar{a}^qTT^*=T^{p+1}{T^*}^{q+1} \in \SC(T, T^*), a^p\bar{a}^qT^*T= {T^*}^{q+1}T^{p+1} \in \SC(T, T^*)$, and for $k>1$: 

$$a^p\bar{a}^q(TT^*)^k = a^p\bar{a}^q(TT^*)(TT^*)^{k-1} = T^{p+1}{T^*}^{q+1}(TT^*)^{k-1} \in \SC(T, T^*)\, \text{and}$$
$$ a^p\bar{a}^q(T^*T)^k = a^p\bar{a}^q(T^*T)(T^*T)^{k-1} = {T^*}^{q+1}T^{p+1}(T^*T)^{k-1} \in \SC(T, T^*).$$
And $\SC(T, T^*)$ being a semigroup, $a^p\bar{a}^q(TT^*)^k T$ and $a^p\bar{a}^q(T^*T)^kT^* \in \SC(T, T^*)$. This verifies the last inclusion and therefore
\begin{equation}\label {E4}
S(T, T^*) = \{a^nT, \bar{a}^nT^*, a^p\bar{a}^q(TT^*)^k, a^p\bar{a}^q(TT^*)^k T, a^p\bar{a}^q(T^*T)^k,  a^p\bar{a}^q(T^*T)^kT^* \mid n\ge0, p, q \ge 0, k\ge 1\}.
\end{equation}

Further simplified form of (\ref{E4}) using Proposition \ref{PT3} becomes 
\begin{equation}\label{E5}
\begin{aligned}
&S(T, T^*)=\\
& \{a^nT, \bar{a}^nT^*, a^p\bar{a}^q||T||^{2k-2}TT^*, a^p\bar{a}^q||T||^{2k}T, a^p\bar{a}^q ||T||^{2k-2}T^*T,  a^p\bar{a}^q||T||^{2k}T^*\mid n\ge0, p, q \ge 0, k\ge 1\}.
\end{aligned}
\end{equation}

These coefficients of the operator products are the only permissible scalars (or any scalars that are representable as one of these), though we are not claiming that as operators they are necessarily different from each other.  
It is also helpful to observe that the monomials in List (\ref{E1})  correspond respectively  with the monomials in Lists (\ref{E4})-(\ref{E5}). \\

(iii) There is a reformulation, which herein we do not explicitly use, of the above simpler forms (i) and (ii) of the semigroup $\SC(T, T^*)$ in terms of projections, with caution to the reader that the projections themselves may not belong to the semigroup but certain scalar multiples of them do. For instance, for any projection $P$, one has $2P \in \SC(2P) = \{(2P)^n, n \ge 1\} = \{2^nP \mid n\geq 1\}$, so $P \not \in \SC(2P)$.

For the list reformulation for rank one operators, recall from Proposition \ref{PT1} that every rank one operator has the form $T = f \otimes T^*f$. It is well known that $P := f \otimes f$ and $Q := \frac{T^*f}{||T^*f||} \otimes \frac{T^*f}{||T^*f||}$ are rank one projections. Indeed, this follows from properties (i)-(ii) of Section 3 (first paragraph), and by simple computations one has $T = PT = TQ$, $TT^* = ||T^*f||^2P = ||T||^2P$ and $T^*T = ||T^*f||^2Q = ||T||^2Q$ and $PQ = \frac{\overline{tr T}}{||T||^2}T$. In the language of projections $P, Q$ together with $T, T^*$, the semigroup lists for $\SC(T, T^*)$ in Remark \ref{RS}(i)-(ii) for $a = 0$ and $a \neq 0$ take the forms:
$$S(T, T^*)= \{0, T, T^*, ||T||^2P,   ||T||^2T, ||T||^2Q,  ||T||^2T^*\}\, \text{and} $$
$$S(T, T^*) = \{a^nT, \bar{a}^nT^*, a^p\bar{a}^q||T||^{2k}P,a^p\bar{a}^q||T||^{2k}T, a^p\bar{a}^q||T||^{2k}Q, a^p\bar{a}^q||T||^{2k}{T^*}, \, n\ge0, p, q \ge 0, k\ge 1\}.$$
Note that $P$ or $Q$
may not belong to $\SC(T, T^*)$, but $PT ( = T)$ and $TQ ( = T^*) \in \SC(T, T^*)$; and scalar multiples of $P$ and $Q$, i.e., $a^p\bar{a}^q||T||^{2k}P$ and $a^p\bar{a}^q||T||^{2k}Q \in \SC(T, T^*)$ suggesting to us that the central role the scalars play here may become a central theme for the general theory of semigroups in $B(H)$ not closed under scalar multiplication.
\eR 

\bR{scalars}
Here we would like to again caution that the role of scalar multiples of operators will be crucial as semigroups may not be closed under scalar multiples which complicates the proof involved in the characterization of selfadjoint-ideal semigroups $\SC(T, T^*)$ generated by rank one $T$ (Theorem \ref{TR2}). For instance, for any  projection $P$, one has $P \in \SC(P) = \{P^n, n \ge 1\} = \{P\}$, but $2P \not \in \SC(P)$.
\eR

\textit{Summary Note --}  The role of Remarks \ref{RS}(ii) and \ref{scalars}: In each case $tr(T)= a =0$ and $ a \neq 0$, there are two forms for the semigroup $\SC(T, T^*)$ - List (\ref{E1}), namely List (\ref{E2}) in terms of $T, T^*, (TT^*)^k, (TT^*)^kT$, $(T^*T)^k$, and  $(T^*T)^kT^*$ for $k \ge 1$ and Lists (\ref{E3})-(\ref{E5}) in terms of scalar multiples of  $T, T^*, (TT^*)^k, (TT^*)^kT$, $(T^*T)^k$, and  $(T^*T)^kT^*$ for $k \ge 1$ and scalar multiples of $T, T^*, TT^*$ and $T^*T$. Which lists we use in our proofs below will depend on the situation.
\vspace{.2cm}

With these simpler forms we are able to characterize the SI property  and simplicity for all semigroups $\SC(T, T^*)$ for rank one operators $T$.
The cases when $tr T = 0$ and $tr T \neq 0$ require separate treatment. The $tr T \neq 0$ case requires by far the most work, and there the non-normal case requires  by far the most work.\\

\noindent \textbf{SI and simple characterizations of $\SC(T, T^*)$ for $T$ rank one: the $tr T =0$ case}\\

For $T$ a rank one operator, the condition $tr T =0$ implies non-normality of $T$. This is because, if $T$ is normal, then
 $T$ is equivalent to 
 \begin{equation*}
 \begin{pmatrix}
a&0\\
0&0
\end{pmatrix} \oplus 0, \text{where $\oplus ~0$ may be absent (proof of Proposition \ref{PN}).}
\end{equation*}
Then $a = tr T \neq 0$, otherwise $T$ is the zero operator against $T$ assumed rank one.

\bT{TR1}
For $T$ a rank one operator with $tr T = 0$, the following are equivalent.
\begin{enumerate}[label=(\roman*)]
\item $\SC(T, T^*)$ is an SI semigroup.
\item $||T|| =1$.
\item $\SC(T, T^*)$ is a simple semigroup.
\end{enumerate}
\eT

\bp
Since $T$ is rank one with $tr T = 0$, the list for $\SC(T, T^*)$ is given by Remark \ref{RS} List (\ref{E3}):  
\begin{equation*}
\SC(T, T^*) = \{0, T, T^*, ||T||^{2n-2}TT^*, ||T||^{2n}T, ||T||^{2n-2}T^*T, ||T||^{2n}T^*\mid  n\ge 1\}.
\end{equation*}
Moreover, $T$ is not normal which follows from the claim above Theorem \ref{TR1}. \\

(i)$\Rightarrow$(ii): If the semigroup $\SC(T, T^*)$ is an SI semigroup, the principal ideal $(T)_{\SC(T, T^*)}$ in particular is selfadjoint and so $T^* \in (T)_{\SC(T, T^*)}$. Because $T$ is not normal, it is not selfadjoint so
$T^* \in (T)_{\SC(T, T^*)}$  implies $T^* = XTY$ where  $X, Y \in \SC(T, T^*)$ and because $T$ is not selfadjoint, at least one $X, Y \neq I$.

\noindent Again as $T$ is not normal, by Lemma \ref{PN}, $Ran T^* \not \subset Ran T$.  So furthermore, because $T^* = XTY$, then $X$ is a word in $T$ and $T^*$ which must start with $T^*$.  
 From these two conditions and because List (\ref{E1}) lists all words in the semigroup $\SC(T, T^*)$, it follows by elimination in List (\ref{E1}) that $XTY$ must be one of the last two words. Then one sees that the first word of these two was reduced in List (\ref{E3}) to a selfadjoint operator, against $T^* = XTY$ not selfadjoint. So  $XTY$ must be the last word of these two, namely, $XTY= (\Pi_{j=1}^{k}{T^*}^{m_j}T^{n_j}){T^*}^{m_{k+1}}$, which in List (\ref{E3}) above reduced to $XTY= ||T||^{2n}T^*$ for some $n \geq 1$. 
Thus $XTY$ has simultaneously the two representations $$T^* = XTY =  ||T||^{2n}T^* \in  (T)_{\SC(T, T^*)}.$$
Since $T^* \neq 0$, so $||T|| =1$.

 (ii)$ \Rightarrow$ (iii): Since $||T|| =1$, the semigroup, again using List (\ref{E3}) above, reduces to $$\SC(T, T^*) = \{0, T, T^*, TT^*, T^*T\}.$$
 
To prove simplicity, it is equivalent to showing that each nonzero principal ideal generated by each operator in $\SC(T, T^*)$ coincides with $\SC(T, T^*)$. Indeed, as $T$ is a partial isomety (Proposition \ref{PT2}), so $T = TT^*T$ and then $T^* = T^*TT^*$. So the principal ideals $(T)_{\SC(T, T^*)}$ and $(T^*)_{\SC(T, T^*)}$ both contain $T$ and $T^*$. Therefore $$(T)_{\SC(T, T^*)} = (T^*)_{\SC(T, T^*)} = \SC(T, T^*)$$ and $T, T^* \in (TT^*)_{\SC(T, T^*)}$ and $T, T^* \in (T^*T)_{\SC(T, T^*)}$ implies $$(TT^*)_{\SC(T, T^*)} = (T^*T)_{\sc(T, T^*)} = \SC(T, T^*).$$

\noindent So the semigroup $\SC(T, T^*)$ contains no proper nonzero ideals, i.e., it is simple.

(iii)$ \Rightarrow $(i)  Simple semigroups contain no nonzero proper ideals, so are automatically SI semigroups.\ep
\vspace{.2cm}

\noindent \textbf{SI and simple characterizations of $\SC(T, T^*)$ for $T$ rank one: the $tr T \neq 0$ case}\\

\textbf{The normal rank one $tr T \neq 0$ case divides into two subcases: selfadjoint and normal nonselfadjoint.}\\
The rank one subcases $T$ selfadjoint and $T$ normal nonselfadjoint require separate treatment in Theorems \ref{Tsa}-\ref{TN1} below. The reason is that the equivalent conditions (i)$\Leftrightarrow$(ii) and (i)$\Leftrightarrow$(iii) in Theorem \ref{TN1} do not hold for $T$ a selfadjoint rank one operator with $tr T \neq 0$. The inequivalence is because: although $\SC(T)$ is automatically SI (Remark \ref{Rsa}(i)), it neither implies that $tr T \in \{-1, 1\}$ nor that $\SC(T)$ is simple (see Remark \ref{R3}(ii)(b)).\\

A normal rank one operator $T$ is equivalent to 
 \begin{equation*}
 \begin{pmatrix}
a&0\\
0&0
\end{pmatrix} \oplus 0, \text{where $\oplus ~0$ may be absent (proof of Proposition \ref{PN}).}
\end{equation*}

For $T$ a \textit{selfadjoint rank one operator with $tr T \neq 0$}, we obtain a characterization of simplicity of the SI semigroup $\SC(T)$ in Theorem \ref{Tsa}. 
\bT{Tsa} 
For $T$ a selfadjoint rank one operator with $tr T \neq 0$, \begin{center}
$\SC(T)$ is simple if and only if $tr T \in \{-1,1\}$. 
\end{center}
\eT
\bp
For $T$ a selfadjoint rank one operator, we use List (\ref{E4}) for $\SC(T, T^*)$ which reduces to $$\SC(T) = \{a^{n}T \mid n\geq 0\}.$$ 
Also, using the above matrix form for $T$ selfadjoint, $0 \neq a = tr T \in \mathbb R$.

Suppose $\SC(T)$ is simple. Then the principal ideal $(aT)_{\SC(T)}$ coincides with $\SC(T)$, i.e.,  $(aT)_{\SC(T)} = \SC(T)$. So $T = X(aT)Y$ where $X, Y \in \SC(T) \cup \{I\}$. Then using the list  $\SC(T) = \{a^{n}T \mid n\geq 0\}$ and the identity $T^n = a^{n-1}T$ (Proposition \ref{PT3}(v)) and $T = aXTY$,
$$T = a^kT$$ for some $k\geq 1$. Since $T \neq 0$, so $a^k =1$, and then since $a \in \mathbb R$, so $a \in \{-1,1\}$.

Conversely, suppose $tr T \in \{-1, 1\}$. For $a = tr T =1$, then using the matrix form of $T$, $\SC(T) = \{T\}$ hence $\SC(T)$ is simple. If $a = tr T = -1$, then using again the matrix form of $T$, $\SC(T) = \{T, T^2\}$ as $T^3 = T$. Since $T = T^2T$, one has $T \in (T^2)_{\SC(T)}$ which implies that $$(T^2)_{\SC(T)} = \{T, T^2\} = \SC(T).$$
Also $T^2 \in (T)_{\SC(T)}$. Since $T, T^2 \in (T)_{\SC(T)}$ so  $$(T)_{\SC(T)} = \SC(T).$$  
Therefore $\SC(T)$ is simple.
\ep

\textit{Normal nonselfadjoint rank one $T$ case.} 
Theorem \ref{N2} applied to a normal nonselfadjoint rank one operator with nonzero trace (i.e., $a \not \in \mathbb R$) simplifies to the following statement in Theorem \ref{TN1} below which we mention without proof as it is a particular case of the more general result proved in Theorem \ref{N2}. (Observe that Condition (ii) of Theorem \ref{N2}, i.e., $T \cong U \oplus 0$, in the rank one case is equivalent to $||T|| = 1 = |tr T|$ which is obvious from the above matrix form.) 

\bT{TN1}
For $T$ a rank one normal nonselfadjoint operator with $tr T \neq 0$, the following are equivalent:
\begin{enumerate}[label=(\roman*)]
\item $\SC(T, T^*)$ is an SI semigroup.
\item $||T|| = |tr T| =1$.
\item $\SC(T, T^*)$ is a simple semigroup.
\end{enumerate}
\eT

%\color{blue}{For a general selfadjoint operator $T$, recall that the SI semigroup $\SC(T)$ may or may not be simple (see examples in Remark \ref{R3}(ii)). 
%}\color{black}
\vspace{.2cm}

\textbf{The non-normal rank one $tr T \neq 0 $ case}\\

As stated earlier, this case requires by far the most work.

In characterizing SI for $\SC(T, T^*)$ for $T$ non-normal rank one $tr T \neq 0 $, we found in ``most" cases that we obtain more, namely, simplicity.  
We first obtain a special class of such SI semigroups $\SC(T, T^*)$ generated by  partial isometries which are not simple (Lemma \ref{rtn} and Remark \ref{RNS} below). All other $\SC(T, T^*)$ in this non-normal rank one $tr T \neq 0 $ case will prove to be simple (they are given by the trace-norm condition case in Theorem \ref{TR2}(ii)). 

In characterizing simplicity versus nonsimplicity for SI semigroups $\SC(T, T^*)$ for $T$ non-normal rank one $tr T \neq 0$, we have the following. Since nonsimplicity implies non-normality (Theorem \ref{N2}), Theorem \ref{Tnonsimple} will be applied directly for this characterization.\\

For $T$ a non-normal rank one $tr T \neq 0$, $T$ is unitarily equivalent to 
 \begin{equation}\label{M}
T = \begin{pmatrix}
a&b\\
0&0
\end{pmatrix} \oplus 0
\end{equation}
where $b \neq 0$. We will refer to this matrix form later in this section.\\

We begin the characterization of this non-normal rank one $tr T \neq 0 $ case with results on partial isometries from this class. From Equation (\ref{M}), this is the class where $a,b \neq 0$.

\bP{Ptr}
For $T\in B(H)$ a non-normal rank one with $tr T \neq 0$, if $T$ is a partial isometry then $|tr T|  \neq 1$. Consequently, a rank one partial isometry with $|tr T |=1$ is normal.
\eP
\bp
Since $T$ is a non-normal rank one with $tr T \neq 0$, $T$ is unitarily equivalent to the matrix form given by (\ref{M}). And since $T$ is a partial isometry, $||T|| =1=\sqrt{|a|^2 + |b|^2}$ with $b\neq 0$. This implies that $|tr T| = |a| \neq 1$. And for $|a| =1$, $b=0$.
\ep

\bL{rtn} 
For $T \in B(H)$ a non-normal rank one partial isometry with $tr T \neq 0$,
$$tr T \in \mathbb R~\text{if and only if}~ \SC(T, T^*)~\text{is an SI semigroup}.$$
\eL
\bp
Since $T$ is a non-normal partial isometry,  $||T|| =1$, $T = TT^*T$ and $T^*T \neq TT^*$. So for $T$ rank one, the semigroup List  (\ref{E5})  of Remark \ref{RS}  for $\SC(T, T^*)$ reduces to the $4$-term list
\begin{equation}\label{EE6}
S(T, T^*) = \{a^p\bar{a}^qTT^*, a^p\bar{a}^qT, a^p\bar{a}^qT^*T, a^p\bar{a}^qT^* \mid p, q \ge 0\}.
\end{equation}
where for convenience we denote $a:= tr T$ and in List  (\ref{E5}) the first two terms are included in the second and fourth terms of List (\ref{EE6}).

$\Rightarrow$: Suppose $tr T \in \mathbb R$. Then semigroup List (\ref{EE6}) for $\SC(T, T^*)$ further simplifies to $$S(T, T^*) = \{a^nT, a^nT^*,  a^nTT^*,  a^nT^*T  \mid n\ge0\}.$$

To show that $\SC(T, T^*)$ is an SI semigroup, it suffices to check that the principal ideals generated by nonselfadjoint operators are selfadjoint (Remark \ref{sa}(ii)). Since $a \in \mathbb R$, the only nonselfadjoint operators in the semigroup $\SC(T, T^*)$ are $\{a^nT, a^nT^* \mid n \ge 0\}$ and hence it suffices to show  $(a^nT)_{\SC(T, T^*)}$  and  $(aT^*)_{\SC(T, T^*)}$ are selfadjoint. 
Recalling the partial isometry defining relation $T^* = T^*TT^*$ one has $$a^nT^* = a^n T^*TT^* = T^*(a^nT)T^* \in (a^nT)_{\SC(T, T^*)},$$ hence the principal ideal $(a^nT)_{\SC(T, T^*)}$ is selfadjoint by Lemma \ref{LP} applied to $\mathcal A = \{a^nT\}$. By Remark \ref{sa}(ii) $2^{\text{nd}}$ paragraph, $(aT)_{\SC(T, T^*)}$ is selfadjoint if and only if $(aT^*)_{\SC(T, T^*)}$ is selfadjoint if and only if $(aT)_{\SC(T, T^*)} = (aT^*)_{\SC(T, T^*)}$. 

$\Leftarrow$: Conversely, suppose $\SC(T, T^*)$ is an SI semigroup. We will use List (\ref{EE6}) above for $\SC(T, T^*)$ to show that $tr T \in \mathbb R$.
Although it may seem natural to consider the selfadjointness of the principal ideal generated by $T$, we found it more useful to instead consider the selfadjointness of the principal ideal generated by $T^2$.
That is, the ideal $(T^2)_{\SC(T, T^*)}$ is selfadjoint. 
So ${T^*}^2 \in (T^2)_{\SC(T, T^*)}$. Therefore ${T^*}^2 = XT^2Y$ where $X, Y \in \SC(T, T^*) \cup \{I\}$ and so by Proposition \ref{PT3}(v), $\bar{a}T^* = aXTY$ and $XTY$ is a scalar multiple of $T^*$. But $T$ being non-normal is equivalent to all nonzero scalar multiples of $T^*$ being non-normal, it follows that $XTY$ is non-normal.  Then since $XTY \in \SC(T, T^*)$ and the first three forms in List (\ref{EE6}) are normal, $XTY$ must be of the last form. Thus $\bar{a}T^* = aXTY = a^p\bar{a}^qT^*$ for some $\mathbf{p\ge 1}$, $q\ge0$. 
As $T^* \neq 0$, 
$$ \bar{a} = a^p\bar{a}^q.$$
The cases $p>1, q \geq 0$ and $p\geq 1, q >0$ do not arise because each of these lead to $p+q \geq 2$ hence $|tr T| =1$. This violates non-normality of $T$ since $T$ is a rank one partial isometry with nonzero trace (Proposition \ref{Ptr}). Therefore, among the cases $p\geq 1, q \geq 0$, only $p=1, q=0$ is an admissible possibility, i.e., $\bar{a} =a$ so $tr T  = a \in \mathbb R$. 
\ep

Although $\SC(T, T^*)$ is SI in Lemma \ref{rtn}, it may not be simple as shown by example in the next remark.

\bR{RNS}(A class of nonsimple SI semigroups $\SC(T, T^*)$ generated by a partial isometry.)
For $T$ as in Lemma \ref{rtn} with $tr T \in \mathbb R$, $\SC(T, T^*)$ is SI but not simple. For $0 <a<1,\, 0 <|b|<1$ and $|a|^2 +|b|^2 = 1$, the $2 \times 2$ matrix
 \begin{equation*}
T = \begin{pmatrix}
a&b\\
0&0
\end{pmatrix}.
\end{equation*}
 is a partial isometry (since $T = TT^*T$, a partial isometry characterization \cite[Corollary 3 of Problem 127]{H}). Since $tr T = a \in \mathbb R$, $\SC(T, T^*)$ is an SI semigroup by Lemma \ref{rtn}. But $\SC(T, T^*)$ is not simple because the principal ideal $(T^2)_{\SC(T, T^*)}$, in particular, is not the entire semigroup $\SC(T, T^*)$. Indeed, assume that $(T^2)_{\SC(T, T^*)} = \SC(T, T^*)$. Then $T = XT^2Y = aXTY$.  Thus $$1 = ||T|| = |a|\,||XTY|| \leq |a|\,||X||\, ||T||\, ||Y|| \leq |a|,$$ hence $|a| \geq 1$ against $0 < a < 1$.
 \eR
 
Lemma \ref{rtn} together with a certain \textit{trace-norm condition} provides next a characterization of SI semigroups $\SC(T, T^*)$ for non-normal rank one trace nonzero operators.

\bT{TR2}

For $T \in B(H)$ a non-normal rank one operator with $tr T \neq 0$, \\
the following are equivalent.
\begin{enumerate}[label=(\roman*)]
\item $\SC(T, T^*)$ is an SI semigroup.
\item Either (but not both) the \textit{trace-norm condition} holds: $$tr^m T\,\, \overline{tr^n T}\, ||T||^{2l} = 1~ \text{for some} ~m,n \ge 0~\text{with}~ m~ \text{or}~n \ge 1~\text{and} ~l \ge 1$$
$$\text{or}$$
$$tr T \in \mathbb R ~\text{and}~ ||T||=1$$ 
\end{enumerate}

\noindent (An equivalent form of the trace-norm condition in (ii) for $m,n \ge 0~\text{with}~ m~ \text{or}~n \ge 1~\text{and} ~l \ge 1$: for $m \neq n$, arg $tr T = \frac{2\pi k}{m-n}$ for some $0 \leq k < |m-n|$ and $|tr T| = ||T||^{\frac{-2l}{m+n}}$; and for $m =n$, $|tr T| = ||T||^{\frac{-l}{m}}$.)
 \eT
\bp
We call the condition in (ii) the \textit{trace-norm condition}: $(tr^mT) (\overline{tr^nT}) ||T||^{2l} = 1$ where $m,n \ge 0$  with $m$ or $n \ge 1$, and $l \ge 1$. We leave the equivalent parenthetical to the reader.
As before, when convenient we use the notation $a := tr T$.\\

(i)$\Rightarrow$(ii): Suppose $\SC(T, T^*)$ is an SI semigroup.
Then the principal ideal $(T)_{\SC(T, T^*)}$, in particular, is selfadjoint so $T^* \in (T)_{\SC(T, T^*)}$. Non-normality of $T$ implies nonselfadjointness, $T^* \neq T$, which, by Lemma \ref{L2} (i)$\Rightarrow$(iii), then implies $T^* = XTY$ where $X$ or $Y$ $\in \SC(T, T^*)\setminus \{I\}$. Since $T$ is non-normal and rank one, by Proposition \ref{PN}, $ran T^* \not \subset ranT$ so $X$ can only be of the form
$$X = {T^*}^p~\text{or}~X = {T^*}^pTX'$$ 
for $p \geq 1$ and $X' \in \SC(T, T^*) \cup \{I\}$.
So $XT$ in $XTY$ is of the form 
$$XT = {T^*}^pT~\text{or}~XT = {T^*}^pTX'T,$$
hence $XTY = {T^*}^pTY$ or $XTY = {T^*}^pTX'TY$, respectively. In summary, $XTY$ at least begins with a $T^*$ and is followed somewhere later by a $T$.
Then using List (\ref{E1}) for $\SC(T, T^*)$, i.e.,
$$
 S(T, T^*) = \{T^n, {T^*}^n,\, \Pi_{j=1}^{k}T^{n_j}{T^*}^{m_j}, \,(\Pi_{j=1}^{k}T^{n_j}{T^*}^{m_j})\,T^{n_{k+1}}, \, \Pi_{j=1}^{k}{T^*}^{m_j}T^{n_j},\, (\Pi_{j=1}^{k}{T^*}^{m_j}T^{n_j})\,{T^*}^{m_{k+1}} \}
 $$
 where $n \ge 1,\,  k\ge1,\, n_j, m_j \ge 1\, \text{for}\, 1 \le j \le k, \,
 \text{and } n_{k+1}, m_{k+1} \geq 1$, one observes that $XTY$ can only correspond to at least one of the fifth or the sixth form in List (\ref{E1}): $$XTY \in \{\ \Pi_{j=1}^{k}{T^*}^{m_j}T^{n_j}, (\Pi_{j=1}^{k}{T^*}^{m_j}T^{n_j}){T^*}^{m_{k+1}}\mid  k\ge1,\, n_j, m_j \ge 1\, \text{for}\, 1 \le j \le k, \, \text{and }m_{k+1} \geq 1\}.$$

Using List (\ref{E5}) in Remark \ref{RS}(ii)  reduced correspondingly from List (\ref{E1}) for $\SC(T, T^*)$: 
\begin{equation*}
\begin{aligned}
&S(T, T^*)=\\
& \{a^nT, \bar{a}^nT^*, a^p\bar{a}^q||T||^{2k-2}TT^*, a^p\bar{a}^q||T||^{2k}T, a^p\bar{a}^q ||T||^{2k-2}T^*T,  a^p\bar{a}^q||T||^{2k}T^*\mid n\ge0, p, q \ge 0, k\ge 1\},
\end{aligned}
\end{equation*}
one extracts the fifth and the sixth terms to obtain $$XTY \in \{a^p\bar{a}^q ||T||^{2k-2}T^*T,  a^p\bar{a}^q||T||^{2k}T^*\mid p, q \ge 0, k\ge 1\}.$$

But $T^* = XTY$ implies $XTY$ has the sixth form since otherwise for the fifth form $T^* = a^p\bar{a}^q||T||^{2k-2}T^*T $ (due to the correspondence between List (\ref{E1}) and List (\ref{E5})), $T$ would be a scalar multiple of a selfadjoint operator, hence normal, against the non-normality of $T$. Therefore, $$T^* = XTY = a^{p}\bar{a}^{q}||T||^{2k}T^*$$ for some fixed $p, q \ge 0, k\ge1$.
Then since $T^* = a^{p}\bar{a}^{q}||T||^{2k}T^*$ and having rank one, $T^* \neq 0$, so $$a^{p}\bar{a}^{q}||T||^{2k}=1.$$
  If either $p$ or $q$ is nonzero, then one obtains the trace-norm condition. In the case $p = q =0$, one has $||T|| =1$ which implies that  $T$ is a partial isometry (Proposition \ref{PT2}). Hence by Lemma \ref{rtn}, $tr T \in \mathbb R$.
Therefore, the SI property for semigroup $\SC(T, T^*)$ implies either the trace-norm condition or $tr T \in \mathbb R$ with $||T|| =1$. To see that not both conditions can be satisfied simultaneously, see Remark \ref{Rtr-norm}.\\

(ii)$\Rightarrow$(i): 
Let $T$ be a rank one, non-normal operator with nonzero trace. We deal with each condition separately and the sufficiency of the trace-norm condition is the hardest proof in this paper.\\

\noindent Case 1. Suppose $0 \neq tr T \in \mathbb R$ and $||T|| =1$. The $||T|| =1$ condition implies $T$ is a partial isometry (Proposition \ref{PT2}) and this together with the real trace condition, by Lemma \ref{rtn}, implies $\SC(T, T^*)$ is an SI semigroup.\\

\noindent Case 2. For the trace-norm condition: $a^m\bar{a}^n||T||^{2l} =1$ for  $m, n \ge 0$, $m$ or $n \ge 1$, and $l \ge 1$, we show  that $\SC(T, T^*)$ is \textbf{simple}, and so vacuously SI.
To show $\SC(T, T^*)$ is simple, we use List (\ref{E5}) in Remark \ref{RS}(ii): 
\begin{align*}
&(\ref{E5})\qquad S(T, T^*)=\\
&\{a^nT, \bar{a}^nT^*, a^p\bar{a}^q||T||^{2k-2}TT^*, a^p\bar{a}^q||T||^{2k}T, a^p\bar{a}^q ||T||^{2k-2}T^*T,  a^p\bar{a}^q||T||^{2k}T^*\mid n\ge0,  p, q \ge 0, k\ge 1\}.
\end{align*}
\vspace{.002cm}

It is important to understand that our approach in what follows reduces noncommutative semigroup operations to somewhat commutative ones in the form of scalar manipulations.\\

To prove simplicity of $\SC(T, T^*)$, clearly it is necessary and sufficient to prove all its principal ideals coincide with $\SC(T, T^*)$. We first prove Claim 1: that the principal ideals generated by $(a^p\bar{a}^q||T||^{2k}T)_{\SC(T, T^*)}$ coincide with $\SC(T, T^*)$ for $p,q \geq 0, k \geq 1$. From this we prove Claim 2: all principal ideals are $\SC(T, T^*)$, and hence simplicity of $\SC(T, T^*)$.\\

\textit{Claim 1.} 
$$(a^p\bar{a}^q||T||^{2k}T)_{\SC(T, T^*)}= \SC(T, T^*)$$
for $p,q \geq 0, k \geq 1$. \\

 \textit{Claim 2.}  
 For each operator $A \in \SC(T, T^*)$, there exists $p,q \geq 0, k \geq 1$ for which $a^p\bar{a}^q||T||^{2k}T \in (A)_{\SC(T, T^*)}$.
 \\
 
\noindent \textit{Claims 1-2} together yield $(A)_{\SC(T, T^*)} = \SC(T, T^*)$ for every nonzero  $A \in \SC(T, T^*)$. Indeed, for each $A \in \SC(T, T^*)$, by Claim 2, $a^p\bar{a}^q||T||^{2k}T \in (A)_{\SC(T, T^*)}$ for some $p, q \geq 0$ and $k\geq 1$ which implies $(a^p\bar{a}^q||T||^{2k}T)_{\SC(T, T^*)} \subset (A)_{\SC(T, T^*)}$. But by Claim 1, one has $(a^p\bar{a}^q||T||^{2k}T)_{\SC(T, T^*)} = \SC(T, T^*)$ for that choice of $p, q$ and $k$. Therefore $(A)_{\SC(T, T^*)} = \SC(T, T^*)$, hence $\SC(T, T^*)$ is simple. Therefore, the trace-norm condition implies the simplicity of $\SC(T, T^*)$. 
\\

\textit{Proof of Claim 1.} Using the trace-norm condition, we show that the principal ideals $$(a^p\bar{a}^q||T||^{2k}T)_{\SC(T, T^*)} = \SC(T, T^*) \qquad \text{for all } p, q \geq 0, k \geq1.$$

To show that $(a^{p}\bar{a}^{q}||T||^{2j}T)_{\SC(T, T^*)} = \SC(T, T^*)$ for $p, q \geq 0$, and $j \geq 1$, it suffices to prove $T, T^* \in (a^{p}\bar{a}^{q}||T||^{2j}T)_{\SC(T, T^*)}$ because then $\SC(T, T^*) \subset (a^{p}\bar{a}^{q}||T||^{2j}T)_{\SC(T, T^*)}$. And together with $(a^{p}\bar{a}^{q}||T||^{2j}T)_{\SC(T, T^*)} \subset \SC(T, T^*)$, one obtains $(a^{p}\bar{a}^{q}||T||^{2j}T)_{\SC(T, T^*)} = \SC(T, T^*)$. \\

Recall here the trace-norm condition $a^m\bar{a}^n||T||^{2l} =1$ for $m, n \ge 0$, $m$ or $n \ge 1$ and $l \ge 1$.

\noindent We analyze the $3$ scalar cases based on $m, n$ in the trace-norm condition:  \\
\noindent Case I. $a^m||T||^{2l} = 1$\, ($n=0, m>0$), \\
Case II. $\bar{a}^n||T||^{2l} = 1$\, ($n>0, m=0$), and \\
Case III. $a^m\bar{a}^n||T||^{2l} = 1$\, ($n, m>0$). \\

\noindent  It suffices to prove only for Cases I and III that $T, T^* \in (a^{p}\bar{a}^{q}||T||^{2j}T)_{\SC(T, T^*)}$ for $p, q \geq 0$ and $j \geq 1$ because Case II is equivalent to Case I since $a^m||T||^{2l} = 1$ if and only if $\bar{a}^m||T||^{2l} = 1$.\\

Case I. $a^m||T||^{2l} =1$ ($n =0, m >0$). Hence $a^{mk}||T||^{2lk} =1$ for all $k \geq 1$. \\

 To prove $T, T^* \in (a^{p}\bar{a}^{q}||T||^{2j}T)_{\SC(T, T^*)}$, 
 
choose $k$ sufficiently large so that $lk-j >0, mk-p-1 >0, mk-q>0$.
$$ \text{Choose}~ \SC(T, T^*) \ni X:= T^*(TT^*)^{lk-j}T^{mk-p} = a^{mk-p-1}||T||^{2lk-2j}T^*T ~ (\text{see Proposition \ref{PT3}(iv)-(v)})~\text{and}$$
$$\SC(T, T^*) \ni Y:= {T^*}^{mk-q+1}(TT^*)^{lk-1} = \bar{a}^{mk-q}||T||^{2lk-2}T^* (\text{again see Proposition \ref{PT3}(iv)-(v)}),$$ one has by Proposition \ref{PT3}(ii),(v): 

\begin{align*}
(a^{p}\bar{a}^{q}||T||^{2j}T)_{\SC(T, T^*)} \ni X a^{p}\bar{a}^{q}||T||^{2j}TY &= a^{mk-p-1}||T||^{2lk-2j}a^{p}\bar{a}^{q}||T||^{2j}\bar{a}^{mk-q}||T||^{2lk-2}T^*T^2T^*\\
&= a^{mk}\bar{a}^{mk} ||T||^{4lk-2}T^*TT^* \qquad (T^2 = aT)\\
&= a^{mk}\bar{a}^{mk}||T||^{4lk}T^* \qquad (T^*TT^* = ||T||^2T^*)\\
&= |a^{mk}||T||^{2lk}|^2T^* = T^* \qquad (a^{mk}||T||^{2lk} =1)
\end{align*}
Therefore $T^* \in (a^{p}\bar{a}^{q}||T||^{2j}T)_{\SC(T, T^*)}$.\\
To obtain this for $T$, use the Case I condition $a^m||T||^{2l} = 1$ and $||T||^{2l}T = (TT^*)^{l}T$ (Proposition \ref{PT3}(iv) here and below) to obtain,
$$T = a^m||T||^{2l}T =a^m (TT^*)(TT^*)^{l-1}T=  (a^mT)T^*(||T||^{2l-2}T).$$
Then since as we have shown $T^* \in (a^{p}\bar{a}^{q}||T||^{2j}T)_{\SC(T, T^*)}$, and  $a^mT = T^{m+1}, ||T||^{2l-2}T = (TT^*)^{l-1}T \in \SC(T, T^*)$, so one obtains $T = (a^mT)T^*(||T||^{2l-2}T) \in (a^{p}\bar{a}^{q}||T||^{2j}T)_{\SC(T, T^*)}$. \\

\vspace{.2cm}

Case II. $\bar{a}^n||T||^{2l} =1$ $(m=0, n>0)$. Recall this is equivalent to Case I.\\

Case III.  $a^m\bar{a}^n||T||^{2l} =1$ $(m,n > 0)$. Hence $a^{mk}\bar{a}^{nk}||T||^{2lk} =1$ if and only if $\bar{a}^{mk}a^{nk}||T||^{2lk} =1$ for $k \geq 1$.\\

To prove $T, T^* \in (a^{p}\bar{a}^{q}||T||^{2j}T)_{\SC(T, T^*)}$,

choose $k$ sufficiently large so that $lk-j >0, mk-p >0, nk-q >0$.

$$ \text{Choose}~ \SC(T, T^*) \ni X:= (TT^*)^{lk-j}T^{m(k+1)-p} = a^{m(k+1)-p-1}||T||^{2lk-2j}T~ (\text{Proposition  \ref{PT3}(iv)-(v)})~ \text{and} $$
$$\SC(T, T^*) \ni Y:= {T^*}^{n(k+1)-q+1}(TT^*)^{l-1}T = \bar{a}^{n(k+1)-q} ||T||^{2l-2}T^*T~ (\text{Proposition  \ref{PT3}(iv)-(v)}),$$
one has by Proposition \ref{PT3}(ii),(v):
\begin{align*}
X(a^{p}\bar{a}^{q}||T||^{2j}T)Y &= a^{m(k+1)-p-1}||T||^{2lk-2j}a^{p}\bar{a}^{q}||T||^{2j}\bar{a}^{n(k+1)-q} ||T||^{2l-2}T^2T^*T\\
&= a^{m(k+1)}\bar{a}^{n(k+1)}||T||^{2lk+2l-2}TT^*T \qquad (T^2 = aT)\\
&= a^{m(k+1)}\bar{a}^{n(k+1)}||T||^{2l(k+1)}T  \qquad (TT^*T = ||T||^2T ~(\text{Proposition \ref{PT3}(iv)}))\\
&= T \qquad (\text{since } a^{m}\bar{a}^{n}||T||^{2l} = 1)
\end{align*}
Therefore $T  = X(a^{p}\bar{a}^{q}||T||^{2j}T)Y \in (a^{p}\bar{a}^{q}||T||^{2j}T)_{\SC(T, T^*)}$.\\ 
To obtain this for $T^*$, using $a^m\bar{a}^n||T||^{2l} =1$ and $||T||^{2l}T^* = (T^*T)^lT^*$ (Proposition \ref{PT3}(ii) here and below), one has 
$$T^* = a^m\bar{a}^n||T||^{2l}T^*= a^m\bar{a}^n(T^*T)^lT^* = a^m\bar{a}^nT^*T(T^*T)^{l-1}T^* = (\bar{a}^nT^*)(a^mT)(||T||^{2l-2}T^*).$$

Since $T \in (a^{p}\bar{a}^{q}||T||^{2j}T)_{\SC(T, T^*)}$ so $a^mT = T^{m+1} \in (a^{p}\bar{a}^{q}||T||^{2j}T)_{\SC(T, T^*)}$. \\
Also  $\bar{a}^nT^* = {T^*}^{n+1}, ~ ||T||^{2l-2}T^* = (T^*T)^{l-1}T^* \in \SC(T, T^*)$ hence $$T^* = (\bar{a}^nT^*)(a^mT)(||T||^{2l-2}T^*) \in (a^{p}\bar{a}^{q}||T||^{2j}T)_{\SC(T, T^*)}.$$ 
 Therefore $T, T^* \in (a^{p}\bar{a}^{q}||T||^{2j}T)_{\SC(T, T^*)}$ implying $(a^{p}\bar{a}^{q}||T||^{2j}T)_{\SC(T, T^*)} = \SC(T, T^*)$.
This completes the proof of \textit{Claim 1}.\\

\noindent \textit{Proof of Claim 2.} 
For $A \in \SC(T, T^*)$ (see List (\ref{E5}) below), using Proposition \ref{PT3}(iv) repeatedly, we prove that for some $p, q \geq 0$, and $k \geq 1$,
\begin{center}
$a^p\bar{a}^q||T||^{2k}T \in (A)_{\SC(T, T^*)}$. 
\end{center}
Recall again
\begin{align*}
&(\ref{E5})\qquad S(T, T^*)=\\
&\{a^nT, \bar{a}^nT^*, a^p\bar{a}^q||T||^{2k-2}TT^*, a^p\bar{a}^q||T||^{2k}T, a^p\bar{a}^q ||T||^{2k-2}T^*T,  a^p\bar{a}^q||T||^{2k}T^*\mid n\ge0, p, q \ge 0, k\ge 1\}.
\end{align*}

\begin{enumerate}[label=(\roman*)]
\item For $A = a^nT$, $(A)_{\SC(T, T^*)} \ni AT^*T =(a^nT)T^*T = a^n||T||^2T $ $(\text{so choose } p=n, q=0, k=1)$.\\
\item For $A= \bar{a}^nT^*$, $(A)_{\SC(T, T^*)} \ni TAT =T(\bar{a}^nT^*)T = \bar{a}^n||T||^2T$ $(\text{so choose } p = 0, q=n, k=1)$.\\
\item For $A = a^p\bar{a}^q||T||^{2k-2}TT^*$, $(A)_{\SC(T, T^*)} \ni AT =  (a^p\bar{a}^q||T||^{2k-2}TT^*)T = a^p\bar{a}^q||T||^{2k}T$.\\

\item For $A = a^p\bar{a}^q ||T||^{2k-2}T^*T$, $(A)_{\SC(T, T^*)}\ni TA = a^p\bar{a}^q ||T||^{2k-2}TT^*T = a^p\bar{a}^q||T||^{2k}T$.\\

\item For $A = a^p\bar{a}^q||T||^{2k}T$, trivially $(A)_{\SC(T, T^*)} \ni A = a^p\bar{a}^q||T||^{2k}T$.\\

\item For $A = a^p\bar{a}^q||T||^{2k}T^*$, $(A)_{\SC(T, T^*)} \ni TAT = T(a^p\bar{a}^q||T||^{2k}T^*)T = a^p\bar{a}^q||T||^{2k+2}T$.
\end{enumerate}
\vspace{.2cm}

This completes the proof of \textit{Claim 2} and the full theorem, except to prove the `or' is exclusive in (ii) of the theorem, see the next remark.
\ep

\bR{Rtr-norm}
For $T$ non-normal rank one $tr T \neq 0$, the \textit{trace-norm condition} and the pair $tr T \in \mathbb R$, $||T|| =1$ are two mutually exclusive conditions. In fact even more, the trace-norm condition and $||T||=1$ conditions are mutually exclusive. Indeed, when $||T|| =1$ and the trace-norm condition is satisfied, then one obtains $|tr T| =1$. That is, $T$ is a rank one partial isometry with $|tr T|=1$ and nonzero trace which implies the normality of $T$, contradicting the assumed non-normality of $T$ (see Proposition \ref{PT2} and \ref{Ptr}).
\eR

\bR{RN}  The trace-norm condition:  $tr^m T\,\, \overline{tr^n T}\, ||T||^{2l} = 1~ \text{for some} ~m,n \ge 0~\text{with} ~m ~\text{or}~n \ge 1~\text{and} ~l \ge 1$ implies $|tr T| \leq 1$ if and only if $||T|| \geq 1$. And recall non-normality of $T$ implies that $|tr T|$ and $||T||$ cannot be simultaneously equal to one, since otherwise $T$ is normal by Remark \ref{Rtr-norm}. So the extreme condition $|tr T| = 1 = ||T||$ occurs only in the normal case. Curiously then even in the rank one normal case,  the trace-norm special case condition is satisfied: $tr T\,\overline{tr T}\,||T||^2 = 1$, albeit this is weaker than the condition $|tr T| = 1 = ||T||$. \eR

Theorem \ref{TR2}  combined with Theorem \ref{TR1} provides an answer to Question \ref{QPI} for partial isometries that are rank one, non-normal and are either a power partial isometry or have a real-valued  trace.

\bC{CPI}
For $T$ a partial isometry that is rank one and non-normal: \begin{enumerate}[label=(\roman*)]
\item If $T$ is a power partial isometry, then $\SC(T, T^*)$ is simple.
\item If $tr T \in \mathbb R$, then $\SC(T, T^*)$ is an SI semigroup.\\
But $\SC(T, T^*)$ can be nonsimple (see Remark \ref{RNS}).
\end{enumerate}
\eC

\bp
(i): Since $T$ is a rank one partial isometry, $||T|| =1$.

Suppose $T$ is also a power partial isometry. Then either $tr T =0$ or $|tr T| = 1$ (Proposition \ref{PI2}). If $tr T = 0$ then  by Theorem \ref{TR1}(ii)$\Rightarrow$(iii), $\SC(T, T^*)$ is simple since $||T|| =1$. On the other hand, if $|tr T| = 1$, then $T$ satisfies the trace-norm condition: $tr T\overline{tr T}||T|| =1$. It then follows from the proof of Theorem \ref{TR2}(ii)$\Rightarrow$(i) that $\SC(T, T^*)$ is simple.

(ii): Suppose $tr T \in \mathbb R$. Since $T$ is a partial isometry, $||T|| =1$. For $T$ a non-normal rank one operator, if $tr T =0$ then Theorem \ref{TR1}(ii)$\Rightarrow$(i) implies that $\SC(T, T^*)$ is an SI semigroup, and if $tr T \neq 0$ then Theorem \ref{TR2}(ii) ($tr T \in \mathbb R$ and $||T|| =1$) implies $\SC(T, T^*)$ is an SI semigroup.
\ep

Considering how often in our study SI semigroups turned out to be simple, as mentioned in the upshot at the end of the introduction, we found it interesting to investigate which SI semigroups are nonsimple. 
For this, concerning the SI semigroups $\SC(T, T^*)$ for nonselfadjoint rank one $T$ we already investigated,
we have the following.\\

\noindent \textbf{Characterization of nonsimple SI semigroup $\SC(T, T^*)$ for nonselfadjoint rank one $T$.}

\noindent The main characterization theorems obtained so far, namely, Theorem \ref{N2}, Theorem \ref{TR1}, and Theorem \ref{TR2} lead us to the following complete characterization of the nonsimple SI semigroups for $\SC(T, T^*)$ generated by a rank one nonselfadjoint operator $T$.
\newpage
\bT{Tnonsimple}
For $T$ a nonselfadjoint rank one operator, the following are equivalent.
\begin{enumerate}[label=(\roman*)]
\item $\SC(T, T^*)$ is a nonsimple SI semigroup.
\item $T$ is non-normal with $tr T \in \mathbb R \setminus \{0, 1, -1\}$ and $||T|| =1$.
\end{enumerate}
\eT
\bp
Suppose $\SC(T, T^*)$ is a nonsimple SI semigroup. Then $T$ must be non-normal. Indeed, if $T$ is normal and $\SC(T, T^*)$ is SI, then $\SC(T, T^*)$ is simple (Theorem \ref{N2}(i)$\Leftrightarrow$(iii)) contradicting the nonsimplicity of $\SC(T, T^*)$. 

To prove that $tr T \in \mathbb R \setminus \{0, 1, -1\}$ and $||T|| =1$, we first show that $tr T \neq 0$. Indeed, if $tr T =0$ then Theorem \ref{TR1}(i)$\Leftrightarrow$(iii) implies that the SI semigroup $\SC(T, T^*)$ is simple, again a contradiction to the nonsimplicity of $\SC(T, T^*)$, so $tr T \neq 0$. 
Since $T$ is rank one, non-normal with $tr T \neq 0$ and $\SC(T, T^*)$ is SI, by Theorem \ref{TR2}(i)$\Rightarrow$(ii), either $T$ satisfies the trace-norm condition or has real trace with $||T|| =1$. But, if $T$ satisfies the trace-norm condition then again by Theorem \ref{TR2} (ii)(trace-norm condition) $\SC(T, T^*)$ turns out to be simple which contradicts our assumption of nonsimplicity of $\SC(T, T^*)$. Therefore, $T$ has real trace with $||T|| =1$.  The condition $||T|| =1$ implies that $T$ is also a partial isometry (Proposition \ref{PT2}). So $T$ being rank one, non-normal, partial isometry with real nonzero trace further implies that $|tr T| \neq 1$ (i.e., $tr T \not \in \{1, -1\}$) (see Proposition \ref{Ptr}). Therefore the nonsimplicity of the SI semigroup $\SC(T, T^*)$ implies $T$ is non-normal with $tr T \in \mathbb R \setminus \{0, 1, -1\}$ and $||T|| =1$.

Conversely, suppose $T$ is non-normal with $tr T \in \mathbb R \setminus \{0, 1, -1\}$ and $||T|| =1$. The conditions $T$ rank one, normal, $tr T \neq 0$, and real trace with $||T|| =1$ implies that $\SC(T, T^*)$ is SI (see Theorem \ref{TR2}(ii)$\Rightarrow$(i)). Finally we prove that $\SC(T, T^*)$ is nonsimple. For this we use List (\ref{E5})  of Remark \ref{RS} for $\SC(T, T^*)$. 
Since $T$ is a rank one non-normal partial isometry with real trace, i.e., $T^*T \neq TT^*$, $||T|| =1$, $T = TT^*T$, and $a :=tr T \in \mathbb R$, so the semigroup List (\ref{E5})  of Remark \ref{RS}  for $\SC(T, T^*)$ reduces to
$$S(T, T^*) = \{a^nT, a^nT^*,  a^nTT^*,  a^nT^*T  \mid n\ge0\}.$$
For nonsimplicity, it suffices to find a principal ideal which is not the full semigroup $\SC(T, T^*)$. We claim that $(aT)_{\SC(T, T^*)}$ is a proper ideal in $\SC(T, T^*)$.  Indeed, suppose otherwise that $(aT)_{\SC(T, T^*)} = \SC(T, T^*)$. Then $T \in (aT)_{\SC(T, T^*)}$. 
By applying Lemma \ref{L1} to $\mathcal A = \{aT\}$ and using the identity $T = TT^*T$, 
 $$(aT)_{\SC(T, T^*)} = \{a^nT, a^nT^*,  a^nTT^*,  a^nT^*T  \mid n\ge 1\}.$$
As $T$ is non-normal and $a \in \mathbb R$, so $T \in (aT)_{\SC(T, T^*)}$ and non-normality of $T$ implies that $T = a^nT$ for some $n \geq1$. Since $a \in \mathbb R$ and $T \neq 0$, one obtains $a^n =1$ and so $a \in  \{ 1, -1\}$ which contradicts our assumption that $a = tr T \in \mathbb R \setminus \{0, 1, -1\}$. Therefore $(aT)_{\SC(T, T^*)}$ is a proper ideal in $\SC(T, T^*)$. \ep

\noindent To reiterate, Theorem \ref{Tnonsimple} characterizes nonsimple SI semigroups $\SC(T, T^*)$ for nonselfadjoint rank one operators $T$.\\

Finally, our characterizations of SI semigroups $\SC(T, T^*)$ for $T$ normal operators and for rank one operators lead further to the classification of these SI semigroups into simple and nonsimple semigroups as follows.

\vspace{.3cm}

\noindent \textbf{Classifications of SI semigroups $\SC(T, T^*)$, simple and nonsimple for $T$ rank one and $T$ normal  with one exception (selfadjoint with rank greater than one case).}\\

\noindent Our classification splits into  cases and subcases:\\

\textit{\noindent T: normal operator}
\begin{enumerate}[label=(\roman*)]
\item $T$ nonselfadjoint: $\SC(T, T^*)$ is simple (hence SI) (Theorem \ref{N2}) 
\item $T$ selfadjoint: SI but may or may not be simple (see Remark \ref{R3}(i)-(ii)  for examples)
\end{enumerate}
\vspace{.2cm}
\newpage
\textit{T: rank one}
\begin{enumerate}[label=(\roman*)]
\item $tr T =0$ case: $\SC(T, T^*)$ is simple (hence SI) (Theorem \ref{TR1}) \\($T$ rank one with zero trace is never normal (\P ~preceding Theorem \ref{TR1}))

\item $tr T \neq 0$ case:
\begin{enumerate}
\item $T$ normal: 
\begin{enumerate}
\item $T$ nonselfadjoint: $\SC(T, T^*)$ is simple (hence SI) (Theorem \ref{N2})
\item $T$ selfadjoint:  $\SC(T, T^*)$ is simple if and only if  $tr T \in \{-1, 1\}$ (Theorem \ref{Tsa})
\end{enumerate}
\item $T$ non-normal:
\begin{enumerate}
\item $\SC(T, T^*)$ is simple iff $T$ satisfies the trace-norm condition (Theorem \ref{TR2})
\item $\SC(T, T^*)$ is nonsimple SI iff $||T|| =1$ with $tr T \in \mathbb R \setminus \{0, 1, -1\}$ (Theorem \ref{Tnonsimple})
\end{enumerate}
\end{enumerate}
\end{enumerate}

\noindent Beyond the rank one case, i.e., for higher rank including infinite rank non-normal operators $T$, we do not have a characterization of SI semigroups $\SC(T, T^*)$.

\bR{RUP}
As mentioned in the Introduction, one upshot of these results is: singly generated SI semigroups $\SC(T, T^*)$ generated by normal and rank one operators $T$ are usually simple and rarely nonsimple. The SI versus simple situation for higher rank operators remains to be studied.
\eR

We end this section using Theorem \ref{TR2} to prove that SI is not similarity invariant, as promised in comments preceding Theorem \ref{T1}.

\bR{example}
While SI is a unitarily invariant property (Theorem \ref{T1}), that SI   is not a similarity invariant can be seen easily using Theorem \ref{TR2}. For the rank one projection,  \begin{equation*}
T = \begin{pmatrix}
1&0\\
0&0
\end{pmatrix},
\end{equation*} because $T = T^*$, $\SC(T, T^*) = \SC(T)$ is an SI semigroup (see paragraph after Definition \ref{D4}). But for the similarity \begin{equation*}
A =\begin{pmatrix}
1&1\\
0&1
\end{pmatrix}T
 \begin{pmatrix}
1 &-1\\
0& 1
\end{pmatrix} = \begin{pmatrix}
1&-1\\
0&0
\end{pmatrix},
\end{equation*} $||A||= \sqrt{2}$ and $tr A = 1$, clearly violating the trace-norm condition of Theorem \ref{TR2}. So $\SC(A , A^*)$ is not an SI semigroup.
\eR

\section{Direct sums of semigroups of $B(H)$, SI and Simplicity}

For semigroups $\SC_1 \subset B(H_1)$, $\SC_2 \subset B(H_2)$,  define the direct sum $\SC_1 \oplus \SC_2:= \{S\oplus T : S \in \SC_1, T \in \SC_2\}  \subset B(H_1 \oplus H_2)$. The direct sum becomes a semigroup with respect to the product defined on $\SC_1 \oplus \SC_2$ as $(X \oplus Y)(Z \oplus W) := XZ \oplus YW$. Moreover, if $\SC_1$ and $\SC_2$ are selfadjoint semigroups, then the adjoint operation on the semigroup $\SC_1 \oplus \SC_2$ defined by $(S_1\oplus S_2)^* := S_1^* \oplus S_2^*$ preserves selfadjointness.

For $\SC_1$ and $\SC_2$ unital semigroups, if $0 \in \SC_2$ then $\SC_1 \oplus \{0\}$ is a nonzero proper ideal of $\SC_1 \oplus \SC_2$; and likewise if $0 \in \SC_1$ then $\{0\} \oplus \SC_2$ is a nonzero proper ideal of $\SC_1 \oplus \SC_2$.  

Our main results in this Section: Direct sums of SI unital semigroups are SI unital semigroups (Proposition \ref{TDS}) and if the semigroups do not contain $0$, then direct sums of simple unital semigroups are simple unital semigroups (Proposition \ref{T3}). \\

\textbf{Inherited simple and unital properties}
  
\bP{T3}\quad 

\indent (i) If $\SC_1, \SC_2$ are unital, then the direct sum $\SC_1 \oplus \SC_2$ is unital. 

(ii)  If $0 \not \in \SC_1, \SC_2$  are simple unital semigroups,  then $\SC_1 \oplus \SC_2$ is a simple unital semigroup. 
\eP
\bp{}
 (i) The semigroup $\SC_1 \oplus \SC_2$ is unital with $I_1 \oplus I_2$ the identity element, where $I_1, I_2$ are identities for $\SC_1$ and $\SC_2$ respectively.

(ii) Suppose $\SC_1, \SC_2$ are simple semigroups. To show $\SC_1 \oplus \SC_2$ is also simple, it suffices to show that if $J$ is an ideal of $\SC_1 \oplus \SC_2$, then the identity element $I_1 \oplus I_2 \in J$.  Let $X_0 \oplus Y_0 \in J$. Since $\SC_1$, $\SC_2$ do not contain $0$, so $X_0, Y_0 \neq 0$. Consider $J_1^{Y_0} := \{X \in \SC_1 \mid X \oplus Y_0 \in J\}$. Since $X_0 \oplus Y_0 \in J$ so $J_1^{Y_0} \neq \emptyset$. We next show that $J_1^{Y_0}$ is an ideal of $\SC_1$. Indeed, for $X \in J^{Y_0}_1$ and $X' \in \SC_1$, $X' \oplus I_2 \in \SC_1 \oplus \SC_2$ and since $J$ is an ideal of $\SC_1 \oplus \SC_2$, $(X' \oplus I_2)(X \oplus Y_0) = X'X \oplus Y_0 \in J$. So $X'X \in J^{Y_0}_1$ and likewise $XX' \in J^{Y_0}_1$, hence $J^{Y_0}_{1}$ is an ideal of $\SC_1$. But $\SC_1$ is a simple semigroup so $J^{Y_0}_1 = \SC_1$ and hence $I_1 \in J^{Y_0}_1$. In particular, $I_1 \oplus Y_0 \in J$. Consider $J_2^{I_1} := \{Y \in \SC_2 \mid I_1 \oplus Y \in J\}$. This set is nonempty as $I_1 \oplus Y_0 \in J$.  Similarly one has $J_2^{I_1} = \SC_2$ so $I_1 \oplus I_2 \in J$. Hence $J = \SC_1 \oplus \SC_2$.
\ep

An easy extension of Proposition \ref{T3}: For a finite collection of simple unital semigroups $\{\SC_1, \cdots, \SC_k\}$ with $0 \not\in \SC_i$ for $1 \leq i \leq k$, $\displaystyle{{\oplus}_{i =1}^{k}} \SC_{i} $ is a simple semigroup.\\

\textbf{Inherited SI property}

\bP{TDS} The unital semigroups $\SC_1, \SC_2$ are SI semigroups if and only if $\SC_1 \oplus \SC_2$ is an SI semigroup.
\eP
\bp
Suppose $\SC_1, \SC_2$ are SI semigroups.  To show that $\SC_1 \oplus \SC_2$ is an SI semigroup, it suffices to show $(S\oplus T)_{\SC_1 \oplus \SC_2}$ is selfadjoint, for every $S \oplus T \in \SC_1 \oplus \SC_2$ (Lemma \ref{L2} (i) $\Leftrightarrow$ (ii)). Since $\SC_1, \SC_2$ are SI semigroups, again by Lemma \ref{L2} (i) $\Leftrightarrow$ (ii) applied to $\SC_1, \SC_2$, the ideals $(S)_{\SC_1}, (T)_{\SC_2}$ are selfadjoint. So $S^* = X_1SY_1$ for $X_1, Y_1 \in \SC_1$ and $T^* = X_2TY_2$ for $X_2, Y_2 \in \SC_2$ (Lemma \ref{L2} (ii)$ \Leftrightarrow $(iii) applied to $\SC_1, \SC_2$). Hence, $$ (X_1 \oplus X_2)(S \oplus T)(Y_1 \oplus Y_2) = X_1SY_1 \oplus X_2TY_2 = S^* \oplus T^*$$
with $X_1 \oplus X_2, Y_1 \oplus Y_2 \in \SC_1 \oplus \SC_2$. Therefore $(S \oplus T)^* \in (S \oplus T)_{\SC_1 \oplus \SC_2}$. 

Conversely, suppose $\SC_1 \oplus \SC_2$ is an SI semigroup. Since $\SC_1, \SC_2$ are unital semigroups, $\SC_1 \oplus \SC_2$ is unital (Proposition \ref{T3}(i)). For $S \in \SC_1$ and $T \in \SC_2$,  $S \oplus T \in \SC_1 \oplus \SC_2$. The ideal $(S\oplus T)_{\SC_1 \oplus \SC_2}$, in particular, is selfadjoint as $\SC_1 \oplus \SC_2$ is SI. Hence, $$S^* \oplus T^* = (S \oplus T)^*= (X_1 \oplus X_2)(S \oplus T)(Y_1 \oplus Y_2) = X_1SY_1 \oplus X_2TY_2$$
with  $X_1 \oplus X_2, Y_1 \oplus Y_2 \in \SC_1 \oplus \SC_2$. Thus $S^* = X_1SY_1$ and $T^* = X_2TY_2$ where $X_1, Y_1 \in \SC_1$ and $X_2, Y_2 \in \SC_2$ implying the principal ideals $(S)_{\SC_1}$ and $(T)_{\SC_2}$ are selfadjoint respectively. Therefore $\SC_1, \SC_2$ are SI semigroups (Lemma \ref{L2} (i) $\Leftrightarrow$ (ii)).
\ep

\end{document}